\documentclass[11pt,letterpaper]{amsart}

\usepackage[dvipsnames]{xcolor}
\usepackage{amsthm,amsfonts,amssymb,amsmath,mathtools,mathrsfs}
\usepackage{enumitem,booktabs}
\usepackage[citecolor=Black,linkcolor=Red,urlcolor=Blue]{hyperref}
\usepackage{cleveref}
\usepackage{tikz}

\newtheorem{theorem}{Theorem}
\newtheorem{proposition}[theorem]{Proposition}
\newtheorem{lemma}[theorem]{Lemma}
\newtheorem{corollary}[theorem]{Corollary}
\newtheorem{theoremintro}{Theorem}

\theoremstyle{definition}
\newtheorem{definition}[theorem]{Definition}
\newtheorem{example}[theorem]{Example}
\newtheorem{remark}[theorem]{Remark}
\numberwithin{equation}{section}
\numberwithin{theorem}{section}

\setlist[itemize]{leftmargin=*,itemsep=\smallskipamount}
\setlist[enumerate]{leftmargin=*,itemsep=\smallskipamount}

\let\ge\geqslant
\newcommand{\kk}{k}
\newcommand{\ZZ}{\mathbb Z}
\newcommand{\QQ}{\mathbb Q}
\newcommand{\CC}{\mathbb C}

\newcommand{\CI}{\mathcal I}

\newcommand{\G}{\mathbb G}
\newcommand{\cH}{\mathcal H}
\newcommand{\bH}{\overline{\mathcal H}}
\newcommand{\cR}{\mathfrak R}
\newcommand{\cK}{\mathfrak K}

\newcommand{\cO}{\mathcal O}
\newcommand{\Wt}{\tilde W}
\newcommand{\Waff}{W_a}
\newcommand{\Saff}{\widetilde{\mathbb S}}
\newcommand{\rig}{\mathrm{rig}}

\newcommand{\naive}{\mathrm{naive}}

\DeclareMathOperator{\Tr}{Tr}
\DeclareMathOperator{\tr}{tr}
\DeclareMathOperator{\Rk}{Rk}

\DeclareMathOperator{\Hom}{Hom}

\DeclareMathOperator{\im}{im}

\DeclareMathOperator{\spanop}{span}

\DeclareMathOperator{\Res}{Res}

\title[Mod-$\ell$ representations of $p$-adic groups]{Modular representations of reductive $p$-adic groups and the kernel of the trace map}

\author[Xuhua He]{Xuhua He}
\address{Department of Mathematics and New Cornerstone Science Laboratory, The University of Hong Kong, Pokfulam, Hong Kong, Hong Kong SAR, China}
\email{xuhuahe@hku.hk}

\newsavebox{\revisionbox}

\DeclareRobustCommand{\revdel}[1]{{%
  \color{red}\hypersetup{allcolors=red}%
  \sbox{\revisionbox}{#1}%
  \leavevmode\rlap{\raisebox{.45ex}{\rule{\wd\revisionbox}{.4pt}}}%
  \usebox{\revisionbox}}}
\hypersetup{colorlinks=true,linkcolor=Black,citecolor=Black}

\begin{document}
\begin{abstract}
The trace map of a finite group gives a perfect duality between conjugacy classes and irreducible complex representations. For reductive $p$-adic groups, an analogous duality holds due to the classical work of Bernstein, Deligne and Kazhdan. In this paper, we study this trace map for modular representations of
reductive $p$-adic groups in characteristic $\ell\ne p$. The failure of the duality is measured by the kernel of the trace map on the cocenter of the ``big Hecke algebra''. 

Our main theorem expresses this kernel in terms of corresponding kernels for finite-dimensional Hecke algebras associated with normalizers of parahoric subgroups in suitable central quotients of Levi subgroups, reducing a problem on infinite-dimensional algebras to finite-dimensional ones. We also show that projective representations provide a natural complementary subspace on the compact part of the cocenter.

For affine Hecke algebras, we prove the naive-kernel conjecture over algebraically closed fields of arbitrary characteristic and for arbitrary compatible nonzero parameters, including unequal parameters and roots of unity.
\end{abstract}



\maketitle

\section*{Introduction}

\subsection{Representations of finite groups}
For a finite group $\Gamma$, the number of irreducible complex
representations equals the number of conjugacy classes.
This equality is a numerical expression of
cocenter--representation duality. The cocenter
\[
 \overline{\CC[\Gamma]}
   =\CC[\Gamma]/[\CC[\Gamma],\CC[\Gamma]]
\]
is the quotient of the group algebra $\CC[\Gamma]$ by its commutator. It has
a standard basis indexed by the conjugacy classes of $\Gamma$.
Let $\cR_\CC(\Gamma)$ be the Grothendieck group of
finite-dimensional complex representations, tensored with
$\CC$. Its standard basis consists of the isomorphism classes of irreducible
representations. The trace map defines a perfect pairing
\[
 \overline{\CC[\Gamma]}\times\cR_\CC(\Gamma)
 \longrightarrow\CC,
 \qquad
 (\bar g,[V])\longmapsto\operatorname{tr}(g\mid V).
\]
With respect to these bases, the matrix of the pairing
is the ordinary character table of $\Gamma$.

In positive characteristic, this duality can fail.
Let $\kk$ be an algebraically closed field of characteristic
$\ell>0$, and define $\cR_\kk(\Gamma)$ analogously.
The trace map
\[
 \Tr_{\kk[\Gamma]}:
 \overline{\kk[\Gamma]}
 \longrightarrow\cR_\kk(\Gamma)^*
\]
remains surjective, but need not be injective.

The cocenter still has a basis indexed by all conjugacy
classes, whereas the irreducible representations are counted
by the $\ell$-regular conjugacy classes, whose elements have order prime to $\ell$.

Let $g_{\ell'}$ be the prime-to-$\ell$ part of $g$.
Every $\kk$-representation has the same trace on $g$ and
$g_{\ell'}$. These relations generate the kernel of the trace map: 
\[
 \ker\Tr_{\kk[\Gamma]}
   =\spanop_\kk
      \{\bar g-\overline{g_{\ell'}}:g\in\Gamma\}.
\]
The $\ell$-regular classes span a complementary subspace
on which the trace map is an isomorphism onto $\cR_\kk(\Gamma)^*$. Thus the number of irreducible $\kk$-representations equals the number of conjugacy classes minus the dimension of the kernel of the trace map. 

\subsection{Cocenter--representation duality for $p$-adic groups}
We study the corresponding questions for a connected reductive $p$-adic group $G$. Its Hecke algebra $H(G)$ plays the role of the group algebra in the finite-group setting. The trace map pairs its cocenter with the Grothendieck group of finite-length smooth representations, tensored with the coefficient field. This space is generally infinite-dimensional.

Over $\CC$, Kazhdan's density theorem asserts that the trace
map is injective, while the trace Paley--Wiener theorem of
Bernstein--Deligne--Kazhdan identifies its image with the
space of good forms \cite{Kaz86,KazClose,BDK}.
Together, these results establish cocenter--representation
duality for reductive $p$-adic groups.

Now let $\kk$ be an algebraically closed field of positive
characteristic $\ell\ne p$. As in the finite-group case,
the trace map need not be injective. Understanding its kernel is therefore essential for extending cocenter--representation duality
to modular representations and for extracting information
about irreducible representations from the cocenter.

\subsection{Main results for $p$-adic groups}
Write $H_\kk=H_\kk(G)$ for the Hecke algebra of $G$.
The compact cocenter $\bar H_\kk^c$ is the image in
$\bar H_\kk$ of the functions supported on compact
elements of $G$.

We first determine which classes in the compact cocenter
have trace zero on every irreducible smooth representation
of $G$. One source of such classes comes from the compact
parts of the normalizers of parahoric subgroups. For each
standard parahoric subgroup $\mathcal P$ and each pro-$p$
level $\CI_n$, the finite-dimensional Hecke algebra
\[
 H_\kk(N_G(\mathcal P)^c,\CI_n)
\]
acts on the space $V^{\CI_n}$ of every irreducible smooth
representation $V$ of $G$. A class in the kernel of its
trace map therefore gives, by extension by zero, a class
in the kernel of $\Tr_\kk^G$. Let $N_\kk(G)$ be the sum
of these images over all $\mathcal P$ and all $n\ge1$.
See \S\ref{sec:kernel} for the notation.

\begin{theoremintro}[{\Cref{cor:compact-kernel-level} and
                      \Cref{thm:rank-regular}}]
\label{thm:intro-compact}
\label{thm:intro-rank}
The kernels of the trace maps for these finite-dimensional
Hecke algebras generate the entire kernel of the trace map
on the compact cocenter:
\[
 \ker\bigl(\Tr_\kk^G|_{\bar H_\kk^c}\bigr)=N_\kk(G).
\]
The compact $\ell$-regular cocenter $\bar H_\kk^{c,\ell'}$
is a natural complement to this kernel. More precisely,
\[
 \bar H_\kk^c=N_\kk(G)\oplus\bar H_\kk^{c,\ell'},
 \qquad
 \bar H_\kk^{c,\ell'}=\im\Rk_\kk,
\]
where $\Rk_\kk$ is the $\kk$-linearized Hattori rank map
for finitely generated projective smooth representations.
\end{theoremintro}

For a finite group, $\overline{g_{\ell'}}$ is the unique
class in the $\ell$-regular cocenter with the same traces
as $\bar g$. The second assertion gives the corresponding
statement for $G$: every compact cocenter class has a unique
compact $\ell$-regular representative with the same traces
on all irreducible smooth representations. Its identification
with the image of the rank map sharpens the abstract Selberg
principle \cite[Theorem~7.1]{CHIII}, which places the image
of the rank map in the compact cocenter. See
\S\ref{sec:rank} for the definitions.

The Newton decomposition \cite{HeI,HeII} writes the cocenter as a direct sum indexed by dominant rational Newton vectors. For each such vector $v$, induction identifies the summand for $G$ with the component
$\bar H_\kk(M_v;v)$ of the rigid cocenter of the corresponding Levi subgroup $M_v$. Our second theorem determines the kernel of the trace map on each summand.

\begin{theoremintro}[{\Cref{thm:whole-kernel}}]
\label{thm:intro-whole}
We have
\[
 \ker\Tr_\kk^G
   =\bigoplus_{v\in V_+}
      \bar i_v\bigl(
        N_\kk^{\rig}(M_v)\cap\bar H_\kk(M_v;v)
      \bigr).
\]
\end{theoremintro}

The condition $x\in N_\kk^{\rig}(M_v)$ in the formula has the following meaning. Push $x\in\bar H_\kk(M_v;v)$ forward to the compact cocenter of a suitable central quotient of $M_v$. Its induced class $\bar i_v(x)$ has trace zero precisely when this pushforward lies in the sum of the kernels contributed by the finite-dimensional Hecke algebras attached to the normalizers of parahoric subgroups of $M_v$. See \S\ref{sec:kernel} for the precise construction.

\subsection{Affine Hecke algebras}

Affine Hecke algebras arise naturally in the representation
theory of reductive $p$-adic groups: for a split reductive $p$-adic group, the convolution algebra of compactly supported, complex-valued
Iwahori-biinvariant functions is an affine Hecke algebra
with equal parameter $q$, where $q$ is the cardinality
of the residue field.
Parameters that are roots of unity are particularly relevant
to modular representations of $p$-adic groups.

For equal complex parameters that are not roots of unity,
Kazhdan and Lusztig \cite{KL} gave a geometric classification
of irreducible representations of affine Hecke algebras.
From the perspective of cocenters, Ciubotaru--He \cite{CH}
established cocenter--representation duality for admissible
complex parameters and trace Paley--Wiener surjectivity
for all nonzero complex parameters.
At roots of unity and in positive characteristic, however,
the representation theory can differ substantially from
the generic complex case. Understanding how the classification of irreducible representations changes in these settings is a major problem in the theory of affine Hecke algebras.

We consider affine Hecke algebras associated with arbitrary
reduced based root data and arbitrary compatible nonzero
parameters, including unequal parameters and roots of unity over an algebraically closed field $\kk$ of arbitrary characteristic.
Write $\cH_\kk^{\mathrm{aff}}$ for the algebra,
$\bH_\kk^{\mathrm{aff}}$ for its cocenter,
and $\Wt$ and $W_0$ for the extended affine Weyl group
and the finite Weyl group, respectively.
The compact cocenter $\bH_\kk^{\mathrm{aff},c}$ has a basis $\{T_{\cO}:\cO\in\operatorname{cl}(\Wt)_0\}$,
where $\operatorname{cl}(\Wt)_0$ is the finite set of
conjugacy classes of finite-order elements of $\Wt$. The indexing set is independent of the parameters.

Let $N_{\naive}(\cH_\kk^{\mathrm{aff}})$ be the sum of
the images in $\bH_\kk^{\mathrm{aff}}$ of the kernels
of the trace maps of the finite extended parahoric
subalgebras $B_{K,\kk}$ defined in \S\ref{sec:affine}.
Our next result proves the naive-kernel conjecture \cite[Conjecture~3.13]{He16}.

\begin{theoremintro}[{\Cref{thm:affine-kernel}}]
\label{thm:intro-affine}
We have $\ker\bigl(
   \Tr_{\cH_\kk^{\mathrm{aff}}}
   |_{\bH_\kk^{\mathrm{aff},c}}
 \bigr)
   =N_{\naive}(\cH_\kk^{\mathrm{aff}})$.
\end{theoremintro}

These results also yield a criterion for injectivity on the whole cocenter of affine Hecke algebras when $|W_0|\in\kk^\times$; see \Cref{thm:inj-affine-whole}.

The theorem also gives a dimension formula for $\cR_\kk(\cH_\kk^{\mathrm{aff}})_{\rig}$, the rigid quotient of the $\kk$-linearized Grothendieck group of finite-dimensional representations.
This quotient identifies parabolically induced classes
that differ by a central twist of the inducing representation;
see \cite[Definition~6.1]{CH}.
Assume that the root datum is semisimple and that
$|W_0|\in\kk^\times$. Then trace induces a perfect pairing
\[
 \bigl(
   \bH_\kk^{\mathrm{aff},c}/
   N_{\naive}(\cH_\kk^{\mathrm{aff}})
 \bigr)
 \times
 \cR_\kk(\cH_\kk^{\mathrm{aff}})_{\rig}
 \longrightarrow\kk.
\]
Consequently, for every compatible family of nonzero
parameters,
\[
 \dim_\kk\cR_\kk(\cH_\kk^{\mathrm{aff}})_{\rig}
   =|\operatorname{cl}(\Wt)_0|
      -\dim_\kk N_{\naive}(\cH_\kk^{\mathrm{aff}}).
\]

For generic complex parameters, the dimension of the rigid
quotient is $|\operatorname{cl}(\Wt)_0|$.
The correction term therefore measures exactly how many
independent irreducible representation classes are
``missing'' in the rigid quotient at roots of unity or
in positive characteristic, in direct analogy with the finite-group counting formula above.

\subsection{Injectivity of the trace map}

The descriptions of the kernel yield explicit criteria for injectivity of the trace map in both settings. 

\begin{theoremintro}[{\Cref{thm:inj-split-whole,cor:affine-injectivity,thm:inj-affine-whole}}]
\label{thm:intro-injectivity}
\leavevmode
\begin{enumerate}
\item
Suppose that $G$ is split and $\ell\ne p$.
The trace map
\[
 \Tr_\kk^G:\bar H_\kk
   \longrightarrow
   \cR_\kk(G)^*_{\mathrm{good}}
\]
is injective if and only if either $\ell=0$, or $\ell$ is banal for $G$.

\item
Suppose that $|W_0|\in\kk^\times$.
For any compatible family of nonzero parameters,
the trace map
\[
 \Tr_{\cH_\kk^{\mathrm{aff}}}:
   \bH_\kk^{\mathrm{aff}}
   \longrightarrow
   \cR_\kk(\cH_\kk^{\mathrm{aff}})^*_{\mathrm{good}}
\]
is injective if and only if every finite parahoric
subalgebra $\cH_{K,\kk}$ is semisimple.
\end{enumerate}
In either case, whenever the trace map is injective,
it is an isomorphism onto the space of good forms.
\end{theoremintro}

A prime $\ell\ne p$ is \emph{banal for $G$} if it does not divide the pro-order of any compact open subgroup of $G$.

For equal parameters and $|W_0|\in\kk^\times$, the criterion for affine Hecke algebras reduces to the condition that $q$ is not a root of the Poincar\'e polynomial of $W_0$.
For equal parameters in $\CC$, this condition is consistent with Xi's result \cite[Theorem~3.3]{Xi} that the Kazhdan--Lusztig classification remains valid whenever
$P_{W_0}(q)\ne0$.

\smallskip

{\bf Acknowledgements} XH is partially supported by the Hong Kong RGC grant 14300023 and by the New Cornerstone Science Foundation through the New Cornerstone Investigator Program. 

\section{Preliminaries}

\label{sec:preliminaries}
\subsection{The \texorpdfstring{$p$}{p}-adic groups and their ``big'' Hecke algebras}
\label{subsec:hecke-algebras}
Let $F$ be a nonarchimedean local field with residue field
$\mathbb F_q$ of characteristic $p$. Let $\G$ be a connected
reductive group over $F$, and put $G=\G(F)$. Let $\kk$ be an
algebraically closed field of characteristic $\ell\ne p$,
where $\ell=0$ is allowed. All representations and modules are over $\kk$.

A representation is smooth if every vector has an open stabilizer.
A smooth representation $V$ is admissible if $V^K$ is
finite-dimensional for every open compact subgroup $K$ of $G$.
Every irreducible smooth representation of $G$ is admissible
by \cite[Proposition~4.10(iv)]{HV}.
Let $\cR(G)$ be the Grothendieck group of finite-length smooth
representations and
$\cR_\kk(G)=\cR(G)\otimes_{\mathbb Z}\kk$.

Fix a maximal $F$-split torus $A$ of $\G$ and an
alcove $\mathfrak a_C$ in its apartment. Let $\CI$ be the
associated Iwahori subgroup and $\CI'$ its pro-$p$ radical.
Choose a minimal parabolic subgroup $P_0=M_0N_0$ with
$M_0=Z_\G(A)(F)$, and put
$W_0=N_\G(A)(F)/Z_\G(A)(F)$.
A parabolic subgroup is called {\it a standard parabolic subgroup} if it contains $P_0$.
The Levi factor containing $M_0$ of a standard parabolic subgroup is called {\it a standard Levi subgroup}.
A parahoric subgroup is called {\it a standard parahoric subgroup} if it contains $\CI$.
We use $P=MN$ for parabolic subgroups and $\mathcal P$
for parahoric subgroups.

The ``big'' Hecke algebra is the convolution algebra
\[
H_\kk=H_\kk(G)=C_c^\infty(G,\kk),\qquad
 (f_1*f_2)(g)=\int_G f_1(x)f_2(x^{-1}g)\, d\mu_G(x).
\]
We normalize the Haar measure $\mu_G$ by $\mu_G(\CI')=1$ and the Haar measure $\mu_M$ on each standard Levi subgroup $M$ by assigning volume $1$ to $M\cap\CI'$.
For an open compact subgroup $K$ of $G$, let $H_\kk(G,K)$ be
the subalgebra of $K$-biinvariant functions. Then
$H_\kk(G)=\varinjlim_K H_\kk(G,K)$, with the subgroups ordered by
reverse inclusion. If $K$ is pro-$p$, put
$e_K=\mu_G(K)^{-1}\mathbf 1_K$. Thus $H_\kk(G,K)=e_K H_\kk(G)e_K$ is
unital with identity $e_K$.

Integration identifies the category of smooth representations
of $G$ with the category of nondegenerate left
$H_\kk(G)$-modules. For an open compact pro-$p$ subgroup $K$,
the functor $V\mapsto e_KV=V^K$ is exact (see \cite[\S3.1]{HV}).
If $K'\subseteq K$ are open compact pro-$p$ subgroups,
then $V^K=e_KV^{K'}$.

The following result is established by Dat--Helm--Kurinczuk--Moss \cite{DHKM}.

\begin{theorem}\label{thm:finiteness}
For every open compact pro-$p$ subgroup $K$ of $G$, the
algebra $H_\kk(G,K)$ is finite over its center, and its center
is a finitely generated $\kk$-algebra.
\end{theorem}

\begin{proof}
For $\ell>0$, this is \cite[Theorem~1.1]{DHKM}.
For $\ell=0$, apply that theorem over $\mathbb Q_{\ell_0}$
for a prime $\ell_0\ne p$. Faithfully flat descent gives
the assertion over $\QQ$ by \cite[Lemma~3.3(2)]{DHKM},
and extension of scalars to $\kk$ gives the result by
\cite[Lemma~3.3(1)]{DHKM}.
\end{proof}

\subsection{The trace map and the trace Paley--Wiener theorem}
\label{subsec:trace}

For a $\kk$-algebra $A$, possibly without identity,
let $[A,A]$ be the $\kk$-linear span of commutators and
put $\bar A=A/[A,A]$. For unital $A$, let $\cR_\kk(A)$ be
the Grothendieck group of finite-dimensional
$A$-modules, tensored with $\kk$. We then have the trace map
$\Tr_A:\bar A\to\cR_\kk(A)^*$.

Put $\bar H_\kk=H_\kk/[H_\kk,H_\kk]$.
Every element of $H_\kk$ acts with finite rank on a finite-length smooth representation. Taking traces gives
\[
 \Tr=\Tr_\kk^G:\bar H_\kk\longrightarrow\cR_\kk(G)^*,
 \qquad \Tr(\bar f,[V])=\tr(f\mid V).
\]
For an open compact pro-$p$ subgroup $K$, the fixed-vector
functor induces a surjection
\[
 \Res_K:\cR_\kk(G)\longrightarrow\cR_\kk(H_\kk(G,K)),
 \qquad [V]\longmapsto[V^K];
\]
see \cite[Theorem~3.2]{HV}.
For $f\in H_\kk(G,K)$, we have
$\tr(f\mid V)=\tr(f\mid V^K)$.
We use the same notation for Levi subgroups.

Fix a square root of $q$ in $\kk$. For a standard parabolic subgroup $P=MN$, let
\[
 i_M^G:\cR_\kk(M)\longrightarrow\cR_\kk(G),\qquad
 r_M^G:\cR_\kk(G)\longrightarrow\cR_\kk(M)
\]
be the normalized parabolic induction and Jacquet restriction
relative to $P$. We suppress $P$ from the notation when
it is understood. The construction in \cite[\S2.2]{CHIII} gives the corresponding maps on the cocenters
\[
 \bar i_M^G:\bar H_\kk(M)\longrightarrow\bar H_\kk(G),\qquad
 \bar r_M^G:\bar H_\kk(G)\longrightarrow\bar H_\kk(M).
\]
Here $\bar r_M^G$ is induced by Jacquet restriction for $P$, and $\bar i_M^G$ is induced by induction from the opposite parabolic. Frobenius reciprocity and second adjointness \cite[Corollary~1.3]{DHKM} give 
\[
 \begin{aligned}
 \Tr_\kk^M(\bar r_M^G(h),\sigma)
   =\Tr_\kk^G(h,i_M^G(\sigma)),\qquad
 \Tr_\kk^G(\bar i_M^G(f),\pi)
   =\Tr_\kk^M(f,r_M^G(\pi)),
 \end{aligned}
\]
for $h\in\bar H_\kk(G)$, $f\in\bar H_\kk(M)$,
$\sigma\in\cR_\kk(M)$ and $\pi\in\cR_\kk(G)$.

Let $M^0$ be the subgroup generated by all compact
subgroups of $M$, and let $\Psi(M)$ be the torus of
characters of the free abelian group $M/M^0$ over $\kk$;
see \cite[\S6.1]{CHIII}.

\begin{definition}\label{def:good-forms}
A linear form $F\in\cR_\kk(G)^*$ is \emph{good} if it
factors through $\Res_K$ for some open compact pro-$p$
subgroup $K$, and if, for every standard parabolic subgroup
$P=MN$ and every $\sigma\in\cR_\kk(M)$, the function
$\chi\mapsto F(i_M^G(\sigma\otimes\chi))$ is regular on
$\Psi(M)$. Denote this space by
$\cR_\kk(G)^*_{\mathrm{good}}$.
\end{definition}

The form $\Tr(h)$ is a good form for any $h \in \bar H_\kk$ (see \cite[\S6.1]{CHIII}). For complex representations, the density and trace Paley--Wiener theorems give the following cocenter--representation duality.

\begin{theorem}\label{thm:classical-trace}
Assume that $\kk=\CC$. Then the trace map induces an isomorphism
\[
 \Tr:\bar H_\kk(G)\xrightarrow{\ \sim\ }
                         \cR_\kk(G)^*_{\mathrm{good}}.
\]
\end{theorem}

Injectivity is Kazhdan's density theorem \cite{Kaz86,KazClose}. For the generality used here, see also \cite[Part~II, \S1.10]{HL}. Surjectivity is the trace Paley--Wiener theorem of Bernstein--Deligne--Kazhdan \cite{BDK}.

The modular trace Paley--Wiener theorem is due to
Ciubotaru and the author \cite[Theorem~6.1]{CHIII}.

\begin{theorem}\label{thm:known-tpw}
Assume that $|W_0|$ is invertible in $\kk$. Then
\[
 \Tr:\bar H_\kk(G)\longrightarrow
                 \cR_\kk(G)^*_{\mathrm{good}}
\]
is surjective.
\end{theorem}

We will use the Newton decomposition to determine $\ker\Tr_\kk^G$ for modular representations. The remaining subsections recall this decomposition and its compatibility with certain levels, then explain how Newton components relate to rigid and compact cocenters.

\subsection{Newton decomposition}
\label{subsec:newton-decomposition}

The Iwahori--Weyl group is
\[
 \Wt=N_\G(A)(F)/(M_0\cap\CI).
\]
The alcove $\mathfrak a_C$ determines
a length function $\ell: \Wt \to \mathbb N$ and a decomposition
$\Wt=\Waff\rtimes\Omega$, where $\Waff$ is the affine
Weyl group and $\Omega \subset \Wt$ is the subgroup of length-zero
elements.
We have the Iwahori--Bruhat decomposition
\[
 G=\bigsqcup_{w\in\Wt}\CI\dot w\CI,
\]
where $\dot w\in N_\G(A)(F)$ represents $w$.

Put $\Lambda=M_0/(M_0\cap\CI)$. This is a finitely generated abelian group that may have torsion. Choosing a special vertex in the closure of $\mathfrak a_C$ gives a decomposition $\Wt\cong\Lambda\rtimes W_0$.
Put $V=\Lambda\otimes_\ZZ\QQ$, naturally
identified with $X_*(A)\otimes_\ZZ\QQ$, and let $V_+$
be its dominant cone determined by $P_0$.

For $w\in\Wt$, its Kottwitz invariant $\kappa(w)\in\Omega$
is its image under the projection $\Wt\to\Omega$.
Its Newton vector records the average translation
of its powers on the apartment. Choose $r\ge1$ such that
$w^r=t^\lambda$ for some $\lambda\in\Lambda$, where $t^\lambda$ denotes the corresponding element of the translation subgroup. The vector $\nu_w=\lambda/r\in V$ is independent of the choice of $r$. Its dominant representative $\bar\nu_w\in V_+$
depends only on the conjugacy class of $w$.

Set $\aleph=\Omega\times V_+$. Following \cite[\S2]{HeI}, we introduce the Newton strata  of  $G$. For $\nu=(\tau,v)\in\aleph$, define the Newton stratum
\[
G(\nu)=
 \bigcup_{\substack{w\in\Wt,\;(\kappa(w),\bar\nu_w)=(\tau,v)\\
                   w\text{ minimal in its conjugacy class}}}
 \ \bigcup_{g\in G}g(\CI\dot w\CI)g^{-1}.
\]
Let $H_\kk(\nu)$ be the space of functions in $H_\kk$
supported on $G(\nu)$, and let $\bar H_\kk(\nu)$ be its
image in $\bar H_\kk$. 

\begin{theorem}[{\cite[Theorems~3 and~10]{HeI}}]
\label{thm:newton-decomposition}
The Newton strata form an open and closed decomposition
$G=\bigsqcup_{\nu\in\aleph}G(\nu)$. We have
\[
H_\kk=\bigoplus_{\nu\in\aleph}H_\kk(\nu),
 \qquad
 \bar H_\kk=\bigoplus_{\nu\in\aleph}\bar H_\kk(\nu).
\]
\end{theorem}

Let $G^c$ be the set of elements whose generated subgroup
has compact closure, and let $G^{\rig}$ be the set of
elements compact modulo $Z(G)$. The images in $\bar H_\kk$
of functions supported on these sets are the compact
cocenter $\bar H_\kk^c$ and the rigid cocenter
$\bar H_\kk^{\rig}$, respectively. Let
$\Omega_{\mathrm{tor}}$ denote the torsion subgroup of
$\Omega$. By \cite[\S6.1]{HeI} and \cite[\S7.1]{CHIII}, we have
\[
 \bar H_\kk^c
   =\bigoplus_{\tau\in\Omega_{\mathrm{tor}}}
      \bar H_\kk((\tau,0)),
 \qquad
 \bar H_\kk^{\rig}
   =\bigoplus_{(\tau,v)\in\aleph, v\in V^{W_0}}
      \bar H_\kk((\tau,v)).
\]

The rigid and compact cocenters also admit descriptions in terms of normalizers of parahoric subgroups.
For a standard parahoric subgroup $\mathcal P$, write $N_G(\mathcal P)$ for its normalizer and $N_G(\mathcal P)^c$ for the subgroup of compact elements of $N_G(\mathcal P)$. This subgroup is the inverse image of the torsion subgroup of $N_G(\mathcal P)/\mathcal P$.
Since this quotient embeds into $\Omega$, $N_G(\mathcal P)^c$ is open compact. The full normalizer $N_G(\mathcal P)$ is compact modulo $Z(G)$. We identify $H_\kk(N_G(\mathcal P)^c)$ and $H_\kk(N_G(\mathcal P))$ as subalgebras of $H_\kk(G)$. By \cite[Proposition~21]{HeI}, we have
\[
 \begin{aligned}
 \bar H_\kk^{\rig}&=\sum_{\mathcal P\supseteq\CI}
       \im\bigl(
         H_\kk(N_G(\mathcal P))
         \to \bar H_\kk
       \bigr),\\
 \bar H_\kk^c
   &=\sum_{\mathcal P\supseteq\CI}
       \im\bigl(
         H_\kk(N_G(\mathcal P)^c)
         \to\bar H_\kk
       \bigr).
 \end{aligned}
\]

\subsection{Newton decomposition at level
\texorpdfstring{$\CI_n$}{I_n}}
\label{subsec:level-newton-decomposition}

For $n\ge1$, let $\CI_n$ be the $n$-th Moy--Prasad
subgroup associated to the barycenter of $\mathfrak a_C$,
as in \cite[\S4.2]{HeI} and \cite[\S1.2]{HeII}.
These subgroups are normal in $\CI$ and form a decreasing
pro-$p$ neighborhood basis of the identity in $G$.
Put $e_n=e_{\CI_n}$ and
$H_{n,\kk}=H_\kk(G,\CI_n)=e_nH_\kk e_n$.
Then $H_{n,\kk}\subseteq H_{m,\kk}$ for $m\ge n$ and
$H_\kk=\bigcup_{n\ge1}H_{n,\kk}$. Set
\[
 \bar H_{n,\kk}
   =\im(H_{n,\kk}\longrightarrow\bar H_\kk)
   \cong H_{n,\kk}/(H_{n,\kk}\cap[H_\kk,H_\kk]).
\]
Then $\bar H_\kk=\bigcup_{n\ge1}\bar H_{n,\kk}$.
The trace map of $H_{n,\kk}$ factors through this image
by the fixed-vector correspondence in
\S\ref{subsec:trace}.

Write $\bar H_{n,\kk}(\nu)$ for the image of
$H_{n,\kk}\cap H_\kk(\nu)$ in $\bar H_\kk$, and define
$\bar H_{n,\kk}^c$ and $\bar H_{n,\kk}^{\rig}$ using
functions in $H_{n,\kk}$ supported on $G^c$ and
$G^{\rig}$, respectively.

Restriction to Newton strata need not preserve
$K$-biinvariance for an open compact subgroup $K$.
Thus the Newton decomposition of $H_\kk$ does not generally restrict to $H_\kk(G,K)$.  Nevertheless, at the levels $\CI_n$, the image $\bar H_{n,\kk}$ in $\bar H_\kk$ admits a Newton decomposition.

\begin{theorem}[{\cite[Theorem~11]{HeI}}]
\label{thm:level-newton-decomposition}
For every $n\ge1$, we have
\[
 \bar H_{n,\kk}
   =\bigoplus_{\nu\in\aleph}\bar H_{n,\kk}(\nu),
 \qquad
 \bar H_{n,\kk}(\nu)=\bar H_{n,\kk}\cap\bar H_\kk(\nu).
\]
Moreover, each $\bar H_{n,\kk}(\nu)$ is finite-dimensional.
\end{theorem}

Finite-dimensionality follows from \cite[Corollary~9 and
Theorem~11(2)]{HeI}. For a fixed $\nu$, only finitely many
minimal-length elements occur, and each corresponding
$\CI\dot w\CI$ contains finitely many $\CI_n$-double cosets.
The decomposition also gives
$\bar H_{n,\kk}^c=\bar H_{n,\kk}\cap\bar H_\kk^c$ and
$\bar H_{n,\kk}^{\rig}=\bar H_{n,\kk}\cap\bar H_\kk^{\rig}$.

By \cite[Theorem~11(2) and Proposition~21]{HeI}, we have
\[
 \begin{aligned}
 \bar H_{n,\kk}^{\rig}
   &=\sum_{\mathcal P\supseteq\CI}
       \im\bigl(
         H_\kk(N_G(\mathcal P),\CI_n)
         \to\bar H_{n,\kk}
       \bigr),\\
 \bar H_{n,\kk}^c
   &=\sum_{\mathcal P\supseteq\CI}
       \im\bigl(
         H_\kk(N_G(\mathcal P)^c,\CI_n)
         \to\bar H_{n,\kk}
       \bigr).
 \end{aligned}
\]

For every standard parahoric subgroup $\mathcal P$,
the algebra $H_\kk(N_G(\mathcal P)^c,\CI_n)$ is
finite-dimensional with identity $e_n$. In particular, $\bar H_{n,\kk}^c$ is finite-dimensional.

For a standard Levi subgroup $M$, put $\CI_M=M\cap\CI$ and
$\CI_{M,n}=M\cap\CI_n$. Then $\CI_M$ is an Iwahori subgroup
of $M$, and $(\CI_{M,n})_{n\ge1}$ is a decreasing pro-$p$
neighborhood basis of the identity in $M$.

\subsection{Reduction to the compact cocenter}
\label{subsec:newton-induction}

Induction identifies each Newton component with a subspace of the rigid cocenter of a standard Levi subgroup. Passing to a discrete central quotient then maps this rigid cocenter to the compact cocenter.

For $v\in V_+$, put
\[
 \bar H_\kk(G;v)
   =\bigoplus_{\tau\in\Omega}\bar H_\kk((\tau,v)),
\]
and use the analogous notation for Levi subgroups, with their own Kottwitz invariants and dominant cones.

For $v\in V_+$, let $M_v$ be the standard Levi
subgroup generated by $M_0$ and the relative root
subgroups for roots $\alpha$ satisfying
$\langle\alpha,v\rangle=0$. Let $P_v=M_vN_v$ be the
standard parabolic subgroup whose unipotent radical
corresponds to the roots with $\langle\alpha,v\rangle>0$.
The vector $v$ is central in $M_v$, so
$\bar H_\kk(M_v;v)\subseteq\bar H_\kk(M_v)^{\rig}$.

\begin{theorem}[{\cite[Theorems~A--C]{HeII}}]
\label{thm:newton-kernel-reduction}
For every $v\in V_+$, the induction map defines an
isomorphism
\[
 \bar i_v:\bar H_\kk(M_v;v)
       \xrightarrow{\ \sim\ }\bar H_\kk(G;v),
\]
adjoint to $r_{M_v}^G$ relative to $P_v$.
Consequently,
\[
 \bar H_\kk
   =\bigoplus_{v\in V_+}
      \bar i_v\bigl(\bar H_\kk(M_v;v)\bigr).
\]
\end{theorem}

We reserve $\bar i_v$ for induction on Newton components
and $\bar i_M^G$ for the categorical induction map on
the full cocenter introduced in \S\ref{subsec:trace}.

Let $M$ be a standard Levi subgroup, including
$M=G$, and let $A_M$ be the maximal $F$-split torus in
its center. Under the identification
$V\simeq X_*(A)\otimes_\ZZ\QQ$, put
$V_{M,c}=X_*(A_M)\otimes_\ZZ\QQ\subseteq V$.
This is the subspace annihilated by all relative roots of $M$. The rigid cocenter therefore decomposes as
\[
 \bar H_\kk(M)^{\rig}
   =\bigoplus_{v\in V_{M,c}}\bar H_\kk(M;v).
\]

Choose a free abelian complement $Z_M$ to the
maximal compact subgroup of $A_M(F)$. For a uniformizer
$\varpi$ of $F$, one may take
$Z_M=\{y(\varpi):y\in X_*(A_M)\}$.
Then $Z_M$ is discrete, and $Z(M)/Z_M$ is compact.
Put $Q_M=M/Z_M$, and let $q_M:M\to Q_M$ be the quotient
map. Equip $Q_M$ with the quotient Haar measure
$\mu_{Q_M}$ characterized by
\[
 \int_{Q_M}\sum_{z\in Z_M}f(mz)\,d\mu_{Q_M}(mZ_M)
   =\int_M f(m)\,d\mu_M(m)
 \qquad(f\in H_\kk(M)).
\]
Pushforward defines a surjective algebra homomorphism
\[
 q_{M,!}:H_\kk(M)\longrightarrow H_\kk(Q_M),
 \qquad
 q_{M,!}(f)(mZ_M)=\sum_{z\in Z_M}f(mz).
\]

Since the map $Q_M\to M/Z(M)$ has compact kernel, an element of $M$ is compact modulo $Z(M)$ if and only if its image in $Q_M$ is compact. Write $\bar q_M:\bar H_\kk(M)\to\bar H_\kk(Q_M)$ for the induced map on cocenters. Its restriction gives a surjection
\[
 \bar q_M:\bar H_\kk(M)^{\rig}
      \longrightarrow\bar H_\kk(Q_M)^c.
\]
For every parahoric subgroup $\mathcal P$ of $M$,
the group $N_M(\mathcal P)/Z_M$ is compact.
The subgroups $q_M(\CI_{M,n})$ are open compact pro-$p$
subgroups of $Q_M$ and form a neighborhood basis of
the identity.

For $v\in V_+$, these constructions give
\[
 \bar H_\kk(G;v)
   \xleftarrow[\sim]{\bar i_v}\bar H_\kk(M_v;v)
   \xrightarrow{\bar q_{M_v}}\bar H_\kk(Q_{M_v})^c,
\]
where $\bar H_\kk(M_v;v)\subseteq\bar H_\kk(M_v)^{\rig}$.
The second arrow is the restriction of central pushforward to this component. 

\section{The kernel of the trace map}

\label{sec:kernel}
\subsection{Description of the kernel}

We describe the kernel of the trace map using finite-dimensional Hecke algebras associated with the normalizers of standard parahoric subgroups in Levi subgroups of $G$, using their compact parts or suitable central quotients.

For a standard parahoric subgroup $\mathcal P$ of $G$,
the algebra $H_\kk(N_G(\mathcal P)^c,\CI_n)$ is
finite-dimensional.
Define the compact naive kernel of $G$ by
\[
 N_\kk(G)=\sum_{\mathcal P,n}
   \im\bigl(
     \ker\Tr_{H_\kk(N_G(\mathcal P)^c,\CI_n)}
     \to\bar H_\kk
   \bigr).
\]

For a standard Levi subgroup $M$, let $Q_M$ be the central quotient constructed in \S\ref{subsec:newton-induction}.
For each standard parahoric subgroup $\mathcal P_M$ of $M$,
the group $N_M(\mathcal P_M)/Z_M$ is compact. Its Hecke
algebra at level $q_M(\CI_{M,n})$ is therefore
finite-dimensional. Define
\[
 N_\kk(Q_M)=\sum_{\mathcal P_M,n}
   \im\bigl(
     \ker\Tr_{H_\kk(N_M(\mathcal P_M)/Z_M,\,q_M(\CI_{M,n}))}
     \to\bar H_\kk(Q_M)
   \bigr).
\]

Using the central quotient map $\bar q_M$, define
\begin{equation}\label{eq:levi-rigid-naive-kernel}
 N_\kk^{\rig}(M)
   =\bigoplus_{v\in V_{M,c}}
      \left\{
        x\in\bar H_\kk(M;v):
        \bar q_M(x)\in N_\kk(Q_M)
      \right\}.
\end{equation}

We call this the rigid naive kernel of $M$.

\begin{theorem}
\label{thm:whole-kernel}
The kernel of the trace map is
\[
 \ker\Tr_\kk^G
   =\bigoplus_{v\in V_+}
      \bar i_v\bigl(
        N_\kk^{\rig}(M_v)\cap\bar H_\kk(M_v;v)
      \bigr).
\]
\end{theorem}

\begin{remark}
\begin{enumerate}
\item For every standard Levi subgroup $M$, the theorem gives
\[
 \ker\!\left(
   \Tr_\kk^{Q_M}\bigm|_{\bar H_\kk(Q_M)^c}
 \right)=N_\kk(Q_M),
 \qquad
 \ker\!\left(
   \Tr_\kk^M\bigm|_{\bar H_\kk(M)^{\rig}}
 \right)=N_\kk^{\rig}(M).
\]
In particular, $N_\kk^{\rig}(M)$ is independent of the
choice of $Z_M$.

\item The sums defining the compact naive kernels need not be direct. Relations among these images are imposed in the ambient cocenter. In contrast,
the direct sum in the theorem is indexed by Newton
vectors, with its $v$-th summand contained in
$\bar H_\kk(G;v)$.
\end{enumerate}
\end{remark}

By \Cref{thm:level-newton-decomposition},
$\bar H_{n,\kk}$ is graded by Newton components.
The formula in \Cref{thm:whole-kernel} therefore determines $\ker(\Tr_\kk^G|_{\bar H_{n,\kk}})$ by intersecting each summand with $\bar H_{n,\kk}$.

\begin{example}\label{ex:gl2-central-quotient}
Let $M=\mathrm{GL}_2(F)$ and
$Z_M=\langle\varpi I_2\rangle$, where $\varpi$ is a
uniformizer of $F$. The element
\[
 w=\begin{pmatrix}0&1\\ \varpi&0\end{pmatrix},
 \qquad w^2=\varpi I_2,
\]
has Newton vector $(\tfrac12,\tfrac12)$.
Central multiplication shifts both slopes by integers,
so no central translate of $w$ is compact.
Its image $c=q_M(w)$ in $Q_M$ is an involution.

The determinant valuation induces a homomorphism
\[
 Q_M\longrightarrow\ZZ/2\ZZ,
 \qquad
 gZ_M\longmapsto\operatorname{val}_F(\det g)\pmod2,
\]
and hence a grading of $H_\kk(Q_M)$ and its cocenter.
For the standard upper triangular Iwahori $\CI$, the
element $w$ normalizes every $\CI_n$.
Let $e_n^Q$ be the normalized idempotent of
$q_M(\CI_n)$, and let $ce_n^Q$ denote its translate by $c$.

In the finite-dimensional algebra $H_\kk(N_M(\CI)/Z_M,q_M(\CI_n))$, we have
$(ce_n^Q-e_n^Q)^2=0$ when $\ell=2$.
This element therefore has trace zero on every
finite-dimensional module over the algebra, and
\[
 \overline{ce_n^Q}-\overline{e_n^Q}\in N_\kk(Q_M).
\]
The two classes have degrees $1$ and $0$, respectively,
and each has trace $1$ on the trivial representation.

Thus projection to a single degree need not preserve $N_\kk(Q_M)$. This explains why \eqref{eq:levi-rigid-naive-kernel} uses the kernel of the trace map on the whole compact cocenter of $Q_M$.
\end{example}

\subsection{Algebraic separation result}
\label{subsec:algebraic-separation}
For a $\kk$-algebra $A$, possibly without identity, put
$E(A)=\spanop_\kk\{\bar e\in\bar A:e\in A,\ e^2=e\}$.
The results in this subsection hold over an algebraically
closed field $\kk$ of arbitrary characteristic.

\begin{lemma}\label{lem:finite}
Let $A$ be a finite-dimensional unital $\kk$-algebra.
The trace map restricts to an isomorphism
\[
 E(A)\xrightarrow{\ \sim\ }\cR_\kk(A)^*.
\]
In particular, $\bar A=\ker\Tr_{A}\oplus E(A)$.
\end{lemma}

\begin{proof}
Let $S_1,\ldots,S_r$ represent the isomorphism classes of
simple $A$-modules. By \cite[Theorem~9.2.1 and its proof]{EtingofEtAl},
we may choose representatives $e_1,\ldots,e_r$ of the conjugacy
classes of primitive idempotents so that $\dim_\kk e_iS_j=\delta_{ij}$.
Every idempotent is a sum of orthogonal primitive idempotents,
and conjugate idempotents have the same image in $\bar A$.
Thus $\bar e_1,\ldots,\bar e_r$ span $E(A)$.
Since $e_i$ acts on $S_j$ as a projection of rank $\delta_{ij}$,
their images under $\Tr_A$ form the basis dual to $[S_1],\ldots,[S_r]$.
\end{proof}

We now extend this separation property to algebras finite over their finitely generated centers.

\begin{theorem}\label{thm:separation}
Let $A$ be a unital $\kk$-algebra and $Z$ be the center of $A$.
Suppose that $Z$ is a finitely generated $\kk$-algebra
and that $A$ is a finite $Z$-module. Then
\[
 E(A)\cap\ker\Tr_A=0.
\]
\end{theorem}

\begin{proof}
Since $z[a,b]=[za,b]$ for $z\in Z$ and $a,b\in A$,
the cocenter $\bar A$ is a finite $Z$-module.
Let $x\in E(A)\cap\ker\Tr_A$, and fix a maximal ideal
$\mathfrak m$ of $Z$. For every $r\ge1$, the algebra
$A/\mathfrak m^rA$ is finite-dimensional over $\kk$:
the Nullstellensatz gives $Z/\mathfrak m=\kk$, and the
successive quotients of the $\mathfrak m$-adic filtration
of $Z/\mathfrak m^r$ are finite-dimensional.
Every simple $A/\mathfrak m^rA$-module is a
finite-dimensional simple $A$-module. Thus the image of $x$ in the cocenter of this quotient has trace zero and is a linear combination of idempotent classes. It therefore vanishes by \Cref{lem:finite}.

The quotient map gives
\[
 \overline{A/\mathfrak m^rA}
   =A/\bigl([A,A]+\mathfrak m^rA\bigr)
   \cong\bar A/\mathfrak m^r\bar A.
\]
Hence $x\in\mathfrak m^r\bar A$ for every $r\ge1$.
The localization $\bar A_{\mathfrak m}
=Z_{\mathfrak m}\otimes_Z\bar A$ is a finite module over
the noetherian local ring $Z_{\mathfrak m}$. By the Krull
intersection theorem,
\[
 x_{\mathfrak m}\in
 \bigcap_{r\ge1}(\mathfrak mZ_{\mathfrak m})^r
                         \bar A_{\mathfrak m}=0.
\]
Since this holds for every maximal ideal $\mathfrak m$,
we obtain $x=0$.
\end{proof}

\begin{remark}\label{rem:thickenings}
The proof uses all infinitesimal quotients $A/\mathfrak m^rA$.
Although their simple modules are the same, their cocenters
may differ.
\end{remark}

For the application to the rank map, we also need a matrix version of the separation property. For a $\kk$-algebra $A$, possibly without identity, the trace of the matrices induces $\operatorname{tr}_r:\overline{M_r(A)}\to\bar A$ for any $r\ge1$. Put
\[
 E_{\mathrm{mat}}(A)
   =\sum_{r\ge1}\operatorname{tr}_r\bigl(E(M_r(A))\bigr).
\]

\begin{corollary}\label{cor:matrix-separation}
Under the hypotheses of \Cref{thm:separation}, we have
$E_{\mathrm{mat}}(A)\cap\ker\Tr_A=0$.
\end{corollary}

\begin{proof}
Let $x \in E_{\mathrm{mat}}(A)\cap\ker\Tr_A$. By adjoining zero blocks, we may write $x=\operatorname{tr}_r(y)$ with $y\in E(M_r(A))$. Every finite-dimensional simple
$M_r(A)$-module is $S^{\oplus r}$ for a finite-dimensional
simple $A$-module $S$. The trace formula for block matrices gives
\[
 \Tr_{M_r(A)}(y,S^{\oplus r})
   =\Tr_A(\operatorname{tr}_r(y),S)
   =\Tr_A(x,S)=0.
\] The algebra $M_r(A)$ is finite over its
center $Z\cdot I_r$. Thus \Cref{thm:separation} gives
$y=0$, and hence $x=0$.
\end{proof}

\subsection{Proof of the formula for the whole cocenter}
\label{subsec:whole-kernel-proof}

We first determine the kernel on the compact cocenter of a central quotient. Fix a standard Levi subgroup $M$, and retain the notation of \S\ref{subsec:newton-induction}. Pushforward induces a surjective algebra homomorphism
\[
 H_\kk(M,\CI_{M,n})
   \longrightarrow H_\kk(Q_M,q_M(\CI_{M,n})).
\]

By \Cref{thm:finiteness}, $H_\kk(M,\CI_{M,n})$ is finite over its center $C$, which is finitely generated. Hence $H_\kk(Q_M,q_M(\CI_{M,n}))$ is finite over $q_{M,!}(C)$, a finitely generated central $\kk$-subalgebra. Since $q_{M,!}(C)$ is noetherian, the
center of $H_\kk(Q_M,q_M(\CI_{M,n}))$ is a finite $q_{M,!}(C)$-module,
and hence is finitely generated over $\kk$. Thus
$H_\kk(Q_M,q_M(\CI_{M,n}))$ satisfies the hypotheses of
\Cref{thm:separation}.

Irreducible smooth representations of $Q_M$ are admissible
by inflation to $M$. The fixed-vector correspondence
\cite[Theorem~3.2]{HV}, which applies
to locally profinite groups, therefore gives a surjection
from $\cR_\kk(Q_M)$ onto
$\cR_\kk(H_\kk(Q_M,q_M(\CI_{M,n})))$.
Write $\Tr_\kk^{Q_M}$ for the resulting trace map on $\bar H_\kk(Q_M)$.

Set $A_n=H_\kk(Q_M,q_M(\CI_{M,n}))$, and let $x\in E(H_\kk(Q_M))\cap\ker\Tr_\kk^{Q_M}$. By definition, we can choose $n\ge1$ and $\widetilde x\in E(A_n)$ in the intrinsic cocenter $A_n/[A_n,A_n]$ so that $x$ is the image of $\widetilde x$ under the map $E(A_n) \to E(H_\kk(Q_M))$.
Surjectivity of the fixed-vector map on Grothendieck
groups and compatibility of traces imply
$\Tr_{A_n}(\widetilde x)=0$.
By \Cref{thm:separation},
$\widetilde x=0$, and hence $x=0$. Thus
\begin{equation}\label{eq:sep-qm}
 E(H_\kk(Q_M))\cap\ker\Tr_\kk^{Q_M}=0.
\end{equation}

Applying $\bar q_M$ to the description of $\bar H_\kk(M)^{\rig}$ gives
\begin{equation}\label{eq:span-qm}
 \bar H_\kk(Q_M)^c
   =\sum_{\mathcal P_M,n}
      \im\bigl(
        H_\kk(N_M(\mathcal P_M)/Z_M,q_M(\CI_{M,n}))
        \longrightarrow\bar H_\kk(Q_M)
      \bigr).
\end{equation}
Here $\mathcal P_M$ ranges over the standard parahoric subgroups of $M$. This uses the surjectivity of $\bar q_M$ on the rigid cocenter and the fact that $H_\kk(N_M(\mathcal P_M)/Z_M)$ is the union of its Hecke algebras at levels $q_M(\CI_{M,n})$.

Restriction to $q_M(\CI_{M,n})$-fixed vectors shows
that $N_\kk(Q_M)\subseteq\ker\Tr_\kk^{Q_M}$.
Conversely, let $x\in\bar H_\kk(Q_M)^c$.
Apply \Cref{lem:finite} to each finite-dimensional algebra
in \eqref{eq:span-qm}. We may then write $x=x_0+y$,
where $x_0\in N_\kk(Q_M)$ and $y\in E(H_\kk(Q_M))$.
If $x$ has trace zero, then so does $y$, and \eqref{eq:sep-qm} implies that $y=0$. Therefore
\[
 \ker\!\left(
   \Tr_\kk^{Q_M}\bigm|_{\bar H_\kk(Q_M)^c}
 \right)=N_\kk(Q_M).
\]

We next use unramified twisting to determine the kernel of the trace map on the rigid cocenter of $M$. Let $\lambda_M:M\to V_{M,c}$ be
the valuation homomorphism characterized by
$\langle\xi,\lambda_M(m)\rangle
=\operatorname{val}_F(\xi(m))$
for every $F$-rational algebraic character $\xi$ of $M$,
with $\operatorname{val}_F(\varpi)=1$.
Its image $\Lambda_M$ is a lattice in $V_{M,c}$.
The restriction of $\lambda_M$ to $Z_M$ is injective,
and $\lambda_M(Z_M)$ has finite index in $\Lambda_M$.
On each Newton stratum with central Newton vector $v$,
we have $\lambda_M=v$. Every
$\chi\in\Hom(\Lambda_M,\kk^\times)$ defines an
unramified character $\chi\circ\lambda_M$ of $M$.

Let $\pi$ be an irreducible smooth representation of $M$,
with central character $\omega_\pi$.
Via $\lambda_M$, its restriction to $Z_M$ defines a character
of the finite-index sublattice $\lambda_M(Z_M)\subseteq\Lambda_M$.
Taking roots in the algebraically closed field $\kk$ extends
this character to $\chi\in\Hom(\Lambda_M,\kk^\times)$.
The twist $\pi\otimes(\chi\circ\lambda_M)^{-1}$ is trivial
on $Z_M$ and hence descends to an irreducible smooth
representation $\rho$ of $Q_M$. Thus $\pi\cong(\rho\circ q_M)\otimes(\chi\circ\lambda_M)$.

Write $x=\sum_v x_v\in\bar H_\kk(M)^{\rig}$, with
$x_v\in\bar H_\kk(M;v)$ and finite support in
$\Lambda_M$. Compatibility of pushforward with
inflation gives
\[
 \Tr_\kk^M\bigl(
   x,(\rho\circ q_M)\otimes(\chi\circ\lambda_M)
 \bigr)
 =\sum_v\chi(v)\,
    \Tr_\kk^{Q_M}\bigl(\bar q_M(x_v),\rho\bigr).
\]
For fixed $\rho$, the functions $\chi\mapsto\chi(v)$ are distinct monomials on the character torus of $\Lambda_M$, hence linearly independent.
Varying $\rho$ and using the preceding description of
irreducible representations and additivity of trace,
we obtain
\[
 x\in\ker\Tr_\kk^M
 \quad\Longleftrightarrow\quad
 \bar q_M(x_v)\in\ker\Tr_\kk^{Q_M}
 \text{ for every }v.
\]
Since each $\bar q_M(x_v)$ belongs to $\bar H_\kk(Q_M)^c$, the formula for the compact cocenter and \eqref{eq:levi-rigid-naive-kernel} give
\[
 \ker\!\left(
   \Tr_\kk^M\bigm|_{\bar H_\kk(M)^{\rig}}
 \right)=N_\kk^{\rig}(M).
\]
This also shows that the rigid naive kernel is independent of the choice of $Z_M$.

We now pass to the whole cocenter. Let $v\in V_+$ and
$x\in N_\kk^{\rig}(M_v)\cap\bar H_\kk(M_v;v)$.
The formula for the rigid cocenter gives $\Tr_\kk^{M_v}(x)=0$.
By adjunction, for every $\pi\in\cR_\kk(G)$,
\[
 \Tr_\kk^G(\bar i_v(x),\pi)
   =\Tr_\kk^{M_v}(x,r_{M_v}^G(\pi))=0.
\]
Thus the right-hand side of \Cref{thm:whole-kernel}
is contained in $\ker\Tr_\kk^G$.

Conversely, let $h\in\ker\Tr_\kk^G$.
By \cite[Theorem~6.3]{HeII}, write
$h=\sum_v h_v$, where
$h_v\in\ker\Tr_\kk^G\cap\bar H_\kk(G;v)$.
Set $x_v=\bar i_v^{-1}(h_v)$, using
\Cref{thm:newton-kernel-reduction}.
By \cite[Corollary~6.2]{HeII},
$x_v\in\ker\Tr_\kk^{M_v}$.
The formula for the rigid cocenter then gives $x_v\in N_\kk^{\rig}(M_v)\cap\bar H_\kk(M_v;v)$.
Hence $h$ belongs to the right-hand side.
The sum is direct by the Newton decomposition.
This proves \Cref{thm:whole-kernel}.

\smallskip

For the compact cocenter, the kernel of the trace map can also be described at each fixed level $\CI_n$.

\begin{corollary}\label{cor:compact-kernel-level}
For every $n\ge1$, we have
\[
 \ker\!\left(
   \Tr_\kk^G\bigm|_{\bar H_{n,\kk}^c}
 \right)
 =\sum_{\mathcal P\supseteq\CI}
   \im\!\left(
     \ker\Tr_{H_\kk(N_G(\mathcal P)^c,\CI_n)}
     \longrightarrow\bar H_{n,\kk}
   \right),
\]
where $\mathcal P$ ranges over the standard parahoric
subgroups of $G$, and the maps are induced by inclusions
of algebras.
\end{corollary}

\begin{proof}
The argument establishing idempotent separation for
$Q_M$ applies to $G$, using \Cref{thm:finiteness}
at the levels $\CI_n$. Hence
$E(H_\kk(G))\cap\ker\Tr_\kk^G=0$.

Restriction to $\CI_n$-fixed vectors shows that the right-hand side is contained in the kernel of the trace map.
Conversely, the compact spanning formula at level
$\CI_n$ in \S\ref{subsec:level-newton-decomposition}
and \Cref{lem:finite} express every
$x\in\bar H_{n,\kk}^c$ as $x=x_0+y$, where $x_0$
belongs to the right-hand side and
$y\in E(H_\kk(G))$.
If $x$ has trace zero, then so does $y$.
Idempotent separation gives $y=0$, proving the assertion.
\end{proof}

\subsection{The rank map on projective modules}\label{sec:rank}

Let $\cK_\kk(G)$ be the Grothendieck group of finitely
generated projective smooth representations, tensored with
$\kk$. Let
\[
 \Rk_\kk:\cK_\kk(G)\longrightarrow\bar H_\kk
\]
be the Hattori rank map, defined as follows. 

Every finitely generated projective smooth representation
$P$ is a direct summand of $(H_\kk e_n)^r$ for some $n,r$.
If $a\in M_r(H_{n,\kk})$ is an idempotent acting on the
right with $P\cong(H_\kk e_n)^{1\times r}a$, then $\Rk_\kk([P])=\overline{\sum_{i=1}^r a_{ii}}$.
Conversely, every idempotent matrix over $H_\kk$ belongs to $M_r(H_{n,\kk})$ for some $n$ and defines such a projective representation. The definition of the rank map therefore gives
\[
 E(H_\kk)\subseteq E_{\mathrm{mat}}(H_\kk)
                    =\im\Rk_\kk,
\]
with the notation of \S\ref{subsec:algebraic-separation}.

For complex representations, the abstract Selberg
principle was proved by Blanc--Brylinski \cite{BB}.
Dat \cite{DatSelberg} gave another proof when
$\operatorname{char}F=0$. The following cocenter
formulation, including the modular case, was established
in the joint work with Ciubotaru
\cite[Theorem~7.1]{CHIII}.

\begin{theorem}
\label{thm:abstract-selberg}
The image of the rank map is contained in the compact
cocenter, i.e.,  $\im\Rk_\kk\subseteq\bar H_\kk^c$.
\end{theorem}

We will prove that $E(H_\kk)=E_{\mathrm{mat}}(H_\kk)=\im\Rk_\kk$ and give an explicit description of this common subspace of $\bar H_\kk^c$.

Let $\Gamma$ be a finite group and $A=\kk[\Gamma]$. Suppose $\ell>0$. For $g\in\Gamma$, write $g=g_\ell g_{\ell'}$ for its
decomposition into commuting $\ell$-power and prime-to-$\ell$
parts. Let
$\Gamma_{\ell'}$ denote the set of elements of order prime to $\ell$. We have the following result. 

\begin{lemma}\label{lem:finite-regular}
Set $\bar A^{\ell'}=\spanop_\kk\{\bar g:g\in\Gamma_{\ell'}\}$. 
We have
\[
 E(A)=E_{\mathrm{mat}}(A)=\bar A^{\ell'},
 \qquad
 \bar A=\ker\Tr_A\oplus\bar A^{\ell'},
\]
where $\ker\Tr_A
   =\spanop_\kk\{\bar g-\overline{g_{\ell'}}:g\in\Gamma\}$. 
\end{lemma}

\begin{proof}
Let $D$ be the span of the differences
$\bar g-\overline{g_{\ell'}}$ for $g\in\Gamma$.
On every finite-dimensional $A$-module, $g_\ell$ acts
unipotently on each eigenspace of $g_{\ell'}$.
Thus $g=g_\ell g_{\ell'}$ and $g_{\ell'}$ have the same
trace, so $D\subseteq\ker\Tr_A$.
Since the classes $\bar g$ span $\bar A$, we have
$\bar A=D+\bar A^{\ell'}$.

For each $\ell$-regular element $g$, the algebra
$\kk[\langle g\rangle]$ is a product of copies of $\kk$.
Hence $g$ is a linear combination of idempotents, and
\[
 \bar A^{\ell'}\subseteq E(A)\subseteq E_{\mathrm{mat}}(A).
\]
By \Cref{cor:matrix-separation}, the trace map is injective
on $E_{\mathrm{mat}}(A)$. It follows that
$\bar A=D\oplus\bar A^{\ell'}$ and $D=\ker\Tr_A$.
Since $E_{\mathrm{mat}}(A)$ contains $\bar A^{\ell'}$ and
meets $D$ trivially, both inclusions above are equalities.
\end{proof}

We now return to the $p$-adic group $G$. 

Suppose $\ell>0$. A compact element $g\in G$ is called
\emph{$\ell$-regular} if every open subgroup of
$\overline{\langle g\rangle}$ has index prime to $\ell$.
Let $G^{c,\ell'}$ be the set of these elements. 

To see that $G^{c, \ell'}$ is open in $G$, recall that every compact element is conjugate into $K=N_G(\mathcal P)^c$ for some standard parahoric subgroup $\mathcal P$.
Choose an open normal pro-$p$ subgroup $J$ of $K$.
Since $\ell\ne p$, an element of $K$ is $\ell$-regular
if and only if its image in $K/J$ has order prime to $\ell$.
Thus $G^{c,\ell'}\cap K$ is a union of open $J$-cosets,
and $G^{c,\ell'}=\cup_{\mathcal P} G \cdot (G^{c,\ell'}\cap N_G(\mathcal P)^c)$ is open. Here $\cdot$ stands for conjugation action.

Let $\bar H_\kk^{c,\ell'}$ be the image in $\bar H_\kk$
of the functions supported on $G^{c,\ell'}$.
For $\ell=0$, set $G^{c,\ell'}=G^c$ and
$\bar H_\kk^{c,\ell'}=\bar H_\kk^c$.

\begin{theorem}
\label{thm:rank-regular}
We have $\im\Rk_\kk=E_{\mathrm{mat}}(H_\kk)=E(H_\kk)=\bar H_\kk^{c,\ell'}$
and $\bar H_\kk^c=N_\kk(G)\oplus E(H_\kk)$.
\end{theorem}

\begin{proof}We first show that
\begin{equation}\label{eq:matrix-trace-separation}
 E_{\mathrm{mat}}(H_\kk)\cap\ker\Tr_\kk^G=\{0\}.
\end{equation}
Let $x$ belong to this intersection. Choosing a common
level for its matrix idempotents gives a lift
$x_n\in E_{\mathrm{mat}}(H_{n,\kk})$ in the intrinsic
cocenter of $H_{n,\kk}$.
Trace compatibility and the surjectivity of
$\Res_{\CI_n}:\cR_\kk(G)\to\cR_\kk(H_{n,\kk})$
give $x_n\in\ker\Tr_{H_{n,\kk}}$.
By \Cref{thm:finiteness} and \Cref{cor:matrix-separation},
we have $x_n=0$, and hence $x=0$.

By the proof of \Cref{cor:compact-kernel-level},
\begin{equation}\label{eq:3-inclusion}
 \bar H_\kk^c
 \subseteq N_\kk(G)+E(H_\kk)
 \subseteq N_\kk(G)+E_{\mathrm{mat}}(H_\kk)
 \subseteq\bar H_\kk^c,
\end{equation}
where the last inclusion follows from the abstract
Selberg principle \Cref{thm:abstract-selberg}.
All inclusions in \eqref{eq:3-inclusion} are therefore equalities.
Since $N_\kk(G)\subseteq\ker\Tr_\kk^G$,
\eqref{eq:matrix-trace-separation} gives
\[
 \bar H_\kk^c=N_\kk(G)\oplus E_{\mathrm{mat}}(H_\kk).
\]
Comparing it with
$\bar H_\kk^c=N_\kk(G)+E(H_\kk)$ and using
$E(H_\kk)\subseteq E_{\mathrm{mat}}(H_\kk)$, we obtain
$E(H_\kk)=E_{\mathrm{mat}}(H_\kk)$.

It remains to identify this common subspace. By the proof of \Cref{cor:compact-kernel-level} and \eqref{eq:matrix-trace-separation},
\[
 E(H_\kk)
   =\sum_{\mathcal P}E\bigl(H_\kk(N_G(\mathcal P)^c)\bigr),
\]
where each term is understood as its image in $\bar H_\kk$.

For any open compact subgroup $K$, passing to finite
quotients with pro-$p$ kernel and applying
\Cref{lem:finite-regular} shows that $E(H_\kk(K))$
is the image of the functions supported on $\ell$-regular
elements. For $\ell=0$, the same assertion holds with
all elements regarded as regular.

Since the conjugates of the $N_G(\mathcal P)^c$ cover
$G^c$ and conjugation acts trivially on the cocenter,
we conclude that $E(H_\kk)=\bar H_\kk^{c,\ell'}$.
\end{proof}

This description also determines the kernel of the rank map. For a finitely generated projective smooth representation $P$ and a finite-length smooth representation $V$, we define a pairing
\[
 \cK_\kk(G)\times\cR_\kk(G)\longrightarrow\kk,
 \qquad
 \langle[P],[V]\rangle
   =\bigl(\dim_\kk\Hom_G(P,V)\bigr)\,1_\kk.
\]
Here the dimension is viewed as an element of $\kk$.

The rank--trace identity \cite[\S1.3]{Dat} gives
\[
 \Tr_\kk^G(\Rk_\kk(x),V)=\langle x,[V]\rangle.
\]
If $\Rk_\kk(x)=0$, then $\langle x,[V]\rangle=0$
for every $V$.
Conversely, if all these pairings vanish, then
$\Rk_\kk(x)$ lies in
$\im\Rk_\kk\cap\ker\Tr_\kk^G$, which is zero by
\Cref{thm:rank-regular}.
Hence $\ker\Rk_\kk$ is precisely the left kernel
of the pairing.

\section{Affine Hecke algebras}

\label{sec:affine}

The arguments of \S\ref{sec:kernel} apply to affine Hecke algebras
associated with arbitrary reduced based root data and
compatible nonzero parameters. Throughout this section, $\kk$ is algebraically closed of arbitrary characteristic $\ell$.

\subsection{The naive kernel}

Let
$\mathscr R=(X,\Phi,X^\vee,\Phi^\vee,\Pi)$ be a reduced
based root datum, and write $\Wt=X\rtimes W_0=\Waff\rtimes\Omega$.
Here $(\Waff,\Saff)$ is the affine Coxeter system and
$\Omega$ is the finitely generated abelian group of
length-zero elements. Write $\Omega_{\mathrm{tor}}$
for its finite torsion subgroup. Let $\pi:\Wt \to W_0$ be the natural
projection, given by $\pi(t^\lambda w)=w$.
Thus $\pi(g)$ is the linear part of the affine
transformation $g\in\Wt$.

Let $\mathbf q=(q_s)_{s\in\Saff}$ be a family of nonzero
parameters in $\kk$ such that $q_s=q_t$ whenever $s$ and $t$
are conjugate in $\Wt$. The affine Hecke algebra
\[
 \cH_\kk^{\mathrm{aff}}
   =\cH_\kk^{\mathrm{aff}}(\Wt,\mathbf q)
\]
has basis $(T_w)_{w\in\Wt}$ and relations
\[
 \begin{aligned}
 T_wT_{w'}&=T_{ww'}
   &&\text{if }\ell(ww')=\ell(w)+\ell(w'),\\
 (T_s-q_s)(T_s+1)&=0
   &&\text{for }s\in\Saff.
 \end{aligned}
\]
Its cocenter is
\[
 \bH_\kk^{\mathrm{aff}}
   =\cH_\kk^{\mathrm{aff}}/
      [\cH_\kk^{\mathrm{aff}},\cH_\kk^{\mathrm{aff}}].
\]

Let $\operatorname{cl}(\Wt)$ denote the set of conjugacy
classes of $\Wt$, and let $\operatorname{cl}(\Wt)_0$
be its finite subset consisting of classes with zero
Newton vector, equivalently, classes of finite-order
elements. For $\cO\in\operatorname{cl}(\Wt)$, let
$T_\cO\in\bH_\kk^{\mathrm{aff}}$ be the class of $T_w$
for a minimal-length element $w\in\cO$.
The combinatorial properties on the minimal-length elements in $\cO$ established in \cite{HN14} show that $T_\cO$ is independent of the choice of $w$. The following theorem gives a basis of the cocenter.

\begin{theorem}
\label{thm:affine-cocenter-basis}
For every compatible family of nonzero parameters,
\[
 \bH_\kk^{\mathrm{aff}}
   =\bigoplus_{\cO\in\operatorname{cl}(\Wt)}
      \kk T_\cO.
\]
\end{theorem}

This basis theorem was proved in the joint work with Nie \cite[Theorem~6.7]{HN14} in the equal-parameter case and in the joint work with Ciubotaru \cite[Theorem~8.1]{CH} for arbitrary nonzero parameters. 

Define the \emph{compact cocenter} by $\bH_\kk^{\mathrm{aff},c}=\bigoplus_{\cO\in\operatorname{cl}(\Wt)_0} \kk T_\cO$. For semisimple root data, this is the rigid cocenter.

A subset $K\subseteq\Saff$ is \emph{spherical} if $W_K$
is finite. For such $K$, let $\cH_{K,\kk}$ be the finite-dimensional Hecke subalgebra generated by $\{T_s:s\in K\}$, and set
\[
 \begin{aligned}
 \Omega_{\mathrm{tor}}(K)
   &=\{\omega\in\Omega_{\mathrm{tor}}:
                        \omega K\omega^{-1}=K\},\\
 B_{K,\kk}
   &=\cH_{K,\kk}\rtimes
                       \kk[\Omega_{\mathrm{tor}}(K)].
 \end{aligned}
\]
Thus $B_{K,\kk}$ is a finite-dimensional unital
subalgebra of $\cH_\kk^{\mathrm{aff}}$. 

These subalgebras play the role of the finite-dimensional Hecke algebras associated with normalizers of parahoric subgroups. For example, let $G$
be a split reductive $p$-adic group with Iwahori--Weyl
group $\Wt$, and let $\CI$ be an Iwahori subgroup.
Take $q_s=q$ for every $s$, where $q$ is the cardinality of the residue field, viewed in $\kk$. Suppose $q, [\CI:\CI']\in\kk^\times$. Then the Iwahori--Matsumoto
presentation identifies
\[
 \cH_\kk^{\mathrm{aff}}\cong H_\kk(G,\CI),
 \qquad
 B_{K,\kk}\cong
 H_\kk(N_G(\mathcal P_K)^c,\CI),
\]
where $\mathcal P_K$ is the standard parahoric subgroup
corresponding to $K$. The identifications send $T_w$ to $\mu_G(\CI)^{-1}\mathbf 1_{\CI\dot w\CI}$, with the Haar normalization of \S\ref{subsec:hecke-algebras}.

Write
$\iota_K:\bar B_{K,\kk}\to\bH_\kk^{\mathrm{aff}}$
for the map induced by inclusion, and set
\[
 \begin{aligned}
 N_{\naive}(\cH_\kk^{\mathrm{aff}})
   &=\sum_{K \text{ is spherical }}\iota_K(\ker\Tr_{B_{K,\kk}}),\\
 E_{\mathrm{par}}(\cH_\kk^{\mathrm{aff}})
   &=\sum_{K \text{ is spherical }}\iota_K(E(B_{K,\kk})).
 \end{aligned}
\]
For semisimple root data, the first equality in the following theorem proves the naive kernel conjecture \cite[Conjecture~3.13]{He16}.

\begin{theorem}\label{thm:affine-kernel}
For every compatible family of nonzero parameters,
\[
 \begin{aligned}
 \ker\bigl(
   \Tr_{\cH_\kk^{\mathrm{aff}}}
      |_{\bH_\kk^{\mathrm{aff},c}}
 \bigr)
   &=N_{\naive}(\cH_\kk^{\mathrm{aff}}),\\
 \bH_\kk^{\mathrm{aff},c}
   &=N_{\naive}(\cH_\kk^{\mathrm{aff}})
      \oplus E_{\mathrm{par}}(\cH_\kk^{\mathrm{aff}}).
 \end{aligned}
\]
\end{theorem}

\begin{proof}
The algebra $\cH_\kk^{\mathrm{aff}}$ is finite over its
finitely generated center; see
\cite[Theorem~2 and \S1.6]{Vig}. By \cite[Corollary~2.8]{HN14} and
\cite[\S3.7.1]{He16}, we have $\bH_\kk^{\mathrm{aff},c}
   =\sum_K\iota_K(\bar B_{K,\kk})$. Restriction to the finite-dimensional subalgebras gives $N_{\naive}(\cH_\kk^{\mathrm{aff}})
   \subseteq\ker\Tr_{\cH_\kk^{\mathrm{aff}}}$.
Applying \Cref{lem:finite} to the finite-dimensional algebras $B_{K,\kk}$, we obtain $\bH_\kk^{\mathrm{aff},c}=N_{\naive}(\cH_\kk^{\mathrm{aff}})+E_{\mathrm{par}}(\cH_\kk^{\mathrm{aff}})$.
Since
$E_{\mathrm{par}}(\cH_\kk^{\mathrm{aff}})
 \subseteq E(\cH_\kk^{\mathrm{aff}})$,
\Cref{thm:separation} shows that the trace map is injective on $E_{\mathrm{par}}(\cH_\kk^{\mathrm{aff}})$. The sum is therefore a direct sum, and its first summand is precisely the kernel of the trace map on the compact cocenter.
\end{proof}

\subsection{The dimension of the rigid quotient}
\label{subsec:affine-characters}

Throughout this subsection, assume that the root datum $\mathscr R$ is semisimple and that $|W_0|\in\kk^\times$. Thus $\bH_\kk^{\mathrm{aff},c}$ is the rigid cocenter. We first recall the rigid quotient introduced in \cite[Definition~6.1]{CH}.

For $J\subseteq\Pi$, let $\cH_{J,\kk}^{\mathrm{aff}}$
be the standard Bernstein parabolic subalgebra and let
$i_J$ denote parabolic induction.

Let
$\cR_\kk(\cH_\kk^{\mathrm{aff}})_{\mathrm{diff}}$
be the subspace of $\cR_\kk(\cH_\kk^{\mathrm{aff}})$ spanned by $i_J(\sigma)-i_J(\sigma\otimes\chi)$ for 
$J\subseteq\Pi$, $\sigma\in\cR_\kk(\cH_{J,\kk}^{\mathrm{aff}})$ and $\chi\in\Hom_{\ZZ}(X/(X\cap\QQ J),\kk^\times)$.
Define the rigid quotient by
\[
 \cR_\kk(\cH_\kk^{\mathrm{aff}})_{\rig}
 =
 \cR_\kk(\cH_\kk^{\mathrm{aff}})
 \big/
 \cR_\kk(\cH_\kk^{\mathrm{aff}})_{\mathrm{diff}}.
\]

\begin{proposition}\label{cor:affine-character-extension}
Assume that the root datum $\mathscr R$ is semisimple
and that $|W_0|\in\kk^\times$. For every compatible family of nonzero parameters, the trace map induces a perfect pairing
\[
 \bigl(\bH_\kk^{\mathrm{aff},c}/
       N_{\naive}(\cH_\kk^{\mathrm{aff}})\bigr)
 \times
 \cR_\kk(\cH_\kk^{\mathrm{aff}})_{\rig}
 \longrightarrow\kk.
\]
In particular, $\dim_\kk\cR_\kk(\cH_\kk^{\mathrm{aff}})_{\rig}
 =
 |\operatorname{cl}(\Wt)_0|
 -
 \dim_\kk N_{\naive}(\cH_\kk^{\mathrm{aff}})$.
\end{proposition}

\begin{proof}
The rigid trace Paley--Wiener theorem (\cite[Theorem~1.1(3)]{CH} and \cite[\S3.7.1]{He16}) gives a surjection
\begin{equation*}
 \Tr_{\cH_\kk^{\mathrm{aff}}}:
 \bH_\kk^{\mathrm{aff},c}
 \twoheadrightarrow
 \cR_\kk(\cH_\kk^{\mathrm{aff}})_{\rig}^{*}.
\end{equation*}

By \Cref{thm:affine-kernel}, the kernel of this restriction of the trace map is
$N_{\naive}(\cH_\kk^{\mathrm{aff}})$. The surjection therefore induces the asserted perfect pairing. The cocenter basis theorem
\Cref{thm:affine-cocenter-basis} gives
$\dim_\kk\bH_\kk^{\mathrm{aff},c}
 =
 |\operatorname{cl}(\Wt)_0|$. This gives the dimension formula.
\end{proof}

The dimension of the rigid quotient is therefore determined by the kernels of the trace maps of the finite Hecke algebras. The examples in \S\ref{sec:examples} illustrate the dimension formula and the relations among the kernels of these finite Hecke algebras.

\par

\section{Injectivity of the trace map}

\label{sec:injectivity}

We use the formulas for the kernel of the trace map on the compact cocenter to characterize injectivity for affine Hecke algebras and split reductive $p$-adic groups. 

\subsection{Injectivity on the compact cocenter $\bH_\kk^{\mathrm{aff},c}$}
\label{subsec:injectivity-affine}

We retain the notation of \S\ref{sec:affine}.

\begin{proposition}\label{cor:affine-injectivity}
The following conditions are equivalent:
\begin{enumerate}
\item The trace map is injective on $\bH_\kk^{\mathrm{aff},c}$.
\item Every $B_{K,\kk}$ is semisimple.
\item $|\Omega_{\mathrm{tor}}|\in\kk^\times$ and
      every $\cH_{K,\kk}$ is semisimple.
\end{enumerate}
\end{proposition}

\begin{proof}
The implication (2)$\Rightarrow$(1) follows from
\Cref{thm:affine-kernel}.
Conversely, the canonical trace
\[
 \tau:
 \text{$\cH_\kk^{\mathrm{aff}}$}
 \longrightarrow\kk,\qquad
 \tau(T_w)=\delta_{w,1},
\]
annihilates commutators and restricts to a symmetrizing
form on each $B_{K,\kk}$.
If $B_{K,\kk}$ is not semisimple, choose a nonzero
element $a$ of its Jacobson radical and
$b\in B_{K,\kk}$ such that $\tau(ab)\ne0$.
Since $ab$ belongs to the radical, it has trace zero
on every finite-dimensional $\cH_\kk^{\mathrm{aff}}$-module.
Its image in $\bH_\kk^{\mathrm{aff},c}$
is nonzero, since it is detected by $\tau$.
This proves (1)$\Rightarrow$(2).

Under (2), the algebra
$B_{\varnothing,\kk}=\kk[\Omega_{\mathrm{tor}}]$
is semisimple, so $|\Omega_{\mathrm{tor}}|\in\kk^\times$.
The radical of $\cH_{K,\kk}$ is
$\Omega_{\mathrm{tor}}(K)$-stable and generates a
nilpotent ideal in $B_{K,\kk}$. Since $B_{K, \kk}$ is semisimple, 
the radical of $\cH_{K,\kk}$ is zero.
Thus (2) implies (3).
Conversely, under (3), averaging module splittings
over $\Omega_{\mathrm{tor}}(K)$ shows that
$B_{K,\kk}$ is semisimple.
\end{proof}

We now consider equal parameters. For a finite Weyl group $W$, let $\cH_W$ be the associated Hecke algebra with parameter $q$ and $P_W(t)=\sum_{w\in W}t^{\ell(w)}$ be the Poincar\'e polynomial. The algebra $\cH_{W,\kk}$ is semisimple if and only if $P_W(q)\ne0$ and either $\ell=0$ or $\ell$ is good for $W$; see
\cite[Proposition~5.3 and the table in \S5.2]{GR}.

For every spherical subset $K\subseteq\Saff$,
the projection map $\pi: \Wt \to W_0$ identifies $W_K$ with a
reflection subgroup of $W_0$. Hence $P_{W_K}(t)\mid P_{W_0}(t)$ in $\ZZ[t]$.
Moreover, a good prime for $W_0$ is good for $W_K$.
The finite-type criterion therefore shows that
semisimplicity of $\cH_{W_0,\kk}$ implies
semisimplicity of every finite parahoric algebra
$\cH_{K,\kk}$.
The converse follows by taking the finite parahoric
algebra $\cH_{W_0,\kk}$. Combining this observation with
\Cref{cor:affine-injectivity} gives the following
criterion.

\begin{corollary}\label{cor:affine-injectivity-equal}
Suppose that $q_s=q$ for all $s\in\Saff$.
The trace map is injective on
$\bH_\kk^{\mathrm{aff},c}$ if and only if $|\Omega_{\mathrm{tor}}|\in\kk^\times$, $P_{W_0}(q)\ne0$, and either $\ell=0$ or $\ell$ is good for $W_0$.
\end{corollary}

\subsection{The whole affine cocenter}
\label{subsec:injectivity-affine-whole}

We retain the notation of
\S\ref{subsec:affine-characters} and \S\ref{subsec:injectivity-affine}.

A linear form
$f\in\cR_\kk(\cH_\kk^{\mathrm{aff}})^*$
is \emph{good} if $\chi\longmapsto f\bigl(i_J(\sigma\otimes\chi)\bigr)$
is regular on
$\Hom_{\ZZ}(X/(X\cap\QQ J),\kk^\times)$
for every $J\subseteq\Pi$ and
$\sigma\in\cR_\kk(\cH_{J,\kk}^{\mathrm{aff}})$.
Write
$\cR_\kk(\cH_\kk^{\mathrm{aff}})^*_{\mathrm{good}}$
for the space of good forms;
cf.\ \cite[Definition~3.5]{CH}.

\begin{theorem}\label{thm:inj-affine-whole}
Assume that $|W_0|\in\kk^\times$ and that
$\cH_{K,\kk}$ is semisimple for every maximal
spherical subset $K\subseteq\Saff$.
Then the trace map induces an isomorphism
\[
 \Tr_{\cH_\kk^{\mathrm{aff}}}:
 \bH_\kk^{\mathrm{aff}}
 \xrightarrow{\;\sim\;}
 \cR_\kk(\cH_\kk^{\mathrm{aff}})^*_{\mathrm{good}}.
\]
In particular, for equal parameters $q_s=q$, the
conclusion holds whenever $|W_0|\in\kk^\times$
and $P_{W_0}(q)\ne0$.
\end{theorem}
\begin{proof}
We prove injectivity by induction on the semisimple
rank, allowing arbitrary root data.
In rank zero, $\cH_\kk^{\mathrm{aff}}$ is the algebra of Laurent
polynomials, and evaluation at characters
separates its elements.
We use throughout that semisimplicity passes from
a finite Hecke algebra to its parabolic subalgebras.

Suppose first that the root datum is semisimple.
The projection $\pi:\Wt\to W_0$ embeds $\Omega$ into $W_0$,
so $|\Omega|\in\kk^\times$.
By \Cref{cor:affine-injectivity}, the trace map is therefore
injective on the compact cocenter.

The hypotheses pass to every proper parabolic subalgebra
$\cH_{J,\kk}^{\mathrm{aff}}$.
We have $|W_J|\in\kk^\times$.
Let $W_{a,J}$ be the affine Coxeter subgroup of
$\Wt_J=X\rtimes W_J$.
Consider a finite parahoric algebra of
$\cH_{J,\kk}^{\mathrm{aff}}$, with finite Coxeter group
$W'\subseteq W_{a,J}$, and choose $x$ whose stabilizer
in $W_{a,J}$ is $W'$.

The ambient stabilizer $W_x=(\Waff)_x$ is a finite
reflection group, acting linearly through $\pi$ when $x$
is taken as the origin. Choose $v$ with stabilizer $W_J$
in $W_0$. The stabilizer $W_x\cap\pi^{-1}(W_J)$ of $v$
in $W_x$ is generated by its reflections fixing $v$.
Their roots belong to the parabolic root subsystem for $J$,
so these affine reflections belong to $W_{a,J}$ and fix $x$.
Thus this stabilizer is contained in $W'$, and the reverse
inclusion follows from the definitions.

Consequently, $W'$ is a parabolic subgroup of $W_x$ and is
conjugate in $\Waff$ to $W_K$ for some spherical
$K\subseteq\Saff$ by \cite[\S1.12 and \S4.8]{Humphreys}.
Conjugation preserves the parameters, so the given finite
parahoric algebra is isomorphic to $\cH_{K,\kk}$ and is
semisimple. The inductive hypothesis now gives injectivity
for $\cH_{J,\kk}^{\mathrm{aff}}$.

Let $h\in\ker\Tr_{\cH_\kk^{\mathrm{aff}}}$.
For every proper subset $J\subsetneq\Pi$, every
$\sigma\in\cR_\kk(\cH_{J,\kk}^{\mathrm{aff}})$,
and every central twist $\chi$, we have
\[
 \Tr_{\cH_\kk^{\mathrm{aff}}}
   \bigl(h,i_J(\sigma\otimes\chi)\bigr)=0.
\]
In particular, these traces are independent of $\chi$.
The separation argument of
\cite[Theorem~7.1]{CH} shows that this independence
forces $h$ to lie in the rigid cocenter:
variation of $\chi$, together with the Mackey
formula, separates the noncentral Newton
contributions, and injectivity for proper parabolic
subalgebras forces these contributions to vanish.
The argument applies over $\kk$, since distinct
Laurent monomials remain linearly independent and
the normalizer quotients used for averaging have
orders dividing $|W_0|$.

Since the root datum is semisimple, the rigid
cocenter equals the compact cocenter.
Injectivity on the compact cocenter therefore
gives $h=0$.

We now consider a general root datum.
Put $L=X^{W_0}$, viewed as the subgroup of central
translations in $\Omega$.
The restriction $\pi|_\Omega$ has kernel $L$.
Since $L$ is torsion-free, we obtain
\begin{equation}\label{eq:inj-affine-central-index}
 \Omega_{\mathrm{tor}}\hookrightarrow
 \Omega/L\hookrightarrow W_0.
\end{equation}
Thus both finite groups have orders invertible
in $\kk$.
The central units $T_z$, $z\in L$, define the quotient
\[
 \cH_\kk^{\mathrm{aff},\mathrm{ss}}
   =\cH_\kk^{\mathrm{aff}}/(T_z-1:z\in L).
\]
Write $\bH_\kk^{\mathrm{aff},\mathrm{ss}}$
for its cocenter.
The quotient has semisimple root datum with lattice
$X/L$, and the same finite Weyl group and finite
parahoric algebras as $\cH_\kk^{\mathrm{aff}}$.
Its trace map is therefore injective by the
semisimple case just proved.

Let $h\in\ker\Tr_{\cH_\kk^{\mathrm{aff}}}$, and write
$h=\sum_{\tau\in\Omega}h_\tau$
for its homogeneous decomposition.
Twisting a representation by a character $\chi$
of $\Omega$ multiplies its trace on $h_\tau$
by $\chi(\tau)$.
Since $|\Omega_{\mathrm{tor}}|\in\kk^\times$,
the functions $\chi\mapsto\chi(\tau)$ are
linearly independent for distinct $\tau$.
It follows that
$\Tr_{\cH_\kk^{\mathrm{aff}}}(h_\tau)=0$
for every $\tau$.

Multiplication by $T_z$ identifies the homogeneous
components whose degrees differ by $z\in L$.
The quotient identifies these components along
their free $L$-orbits, so the induced map
\begin{equation}\label{eq:inj-affine-degreewise}
 \bH_\kk^{\mathrm{aff}}
 \longrightarrow
 \bH_\kk^{\mathrm{aff},\mathrm{ss}}
\end{equation}
is injective on each homogeneous component.
Inflating representations of
$\cH_\kk^{\mathrm{aff},\mathrm{ss}}$
shows that the image of $h_\tau$ has trace zero.
Injectivity of the trace map for the quotient
algebra therefore makes this image zero.
By \eqref{eq:inj-affine-degreewise},
$h_\tau=0$ for every $\tau$, and hence $h=0$.
This proves injectivity.

Surjectivity follows from the argument of \cite[Theorem~8.5]{CH}.
As in the $p$-adic case of \cite[\S6]{CHIII}, the hypothesis
$|W_0|\in\kk^\times$ ensures that the normalizing factors in
the induction--restriction projectors are invertible in $\kk$.
Finally, $|W_0|\in\kk^\times$ implies that $\ell=0$
or $\ell$ is good for $W_0$, so the equal-parameter
assertion follows from
\Cref{cor:affine-injectivity-equal}.
\end{proof}

Over $\CC$, the joint work with Ciubotaru
\cite[Definition~6.6 and Theorem~8.2]{CH}
proved injectivity for admissible parameters.
For semisimple root data and equal complex parameters,
Dawydiak \cite[Theorem~9]{DawRig} proved that the rigid
pairing is perfect whenever $P_{W_0}(q)\ne0$.
Theorem~\ref{thm:inj-affine-whole} gives a criterion
in terms of finite parahoric algebras for unequal parameters and extends injectivity to positive characteristic under the stated hypotheses.

A related result concerns Lusztig's map to the
asymptotic algebra.
For equal complex parameters, this map induces an
isomorphism on cocenters whenever $q$ is not a root
of unity \cite[Theorem~1 and Appendix]{BDD23}.
It remains surjective on cocenters whenever
$P_{W_0}(q)\ne0$ \cite[Remark~2]{BDD23}.

In positive characteristic, the criterion simplifies
when all parameters belong to $\mathbb F_\ell^\times$.
In this case, semisimplicity of $\cH_{W_0,\kk}$
already implies $|W_0|\in\kk^\times$, by
\cite[Proposition~5.3 and the table in \S5.2]{GR}.
Thus \Cref{cor:affine-injectivity,thm:inj-affine-whole}
show that the trace map is an isomorphism onto the
good forms if and only if every finite parahoric
algebra $\cH_{K,\kk}$ is semisimple.
For equal parameters, this is equivalent to
$P_{W_0}(q)\ne0$.
This specialization applies to the Iwahori Hecke
algebras of split reductive $p$-adic groups:
their parameter is the cardinality of the residue field,
whose image lies in $\mathbb F_\ell^\times$
when $\ell\ne p$.

\subsection{Split reductive \texorpdfstring{$p$}{p}-adic groups}
\label{subsec:injectivity-split}

Let $G$ be the group of $F$-points of a split
connected reductive group, and assume that $\ell\ne p$.
For $\ell>0$, recall that $\ell$ is \emph{banal for $G$}
if it does not divide the pro-order of any compact
open subgroup of $G$.
Equivalently, $\ell\nmid |K/U|$ for every compact
open subgroup $K$ and every open normal pro-$p$
subgroup $U$ of $K$.

\begin{theorem}\label{thm:inj-split-whole}
The trace map $\Tr_\kk^G$ on $\bar H_\kk(G)$ is
injective if and only if either $\ell=0$, or
$\ell>0$ and $\ell$ is banal for $G$.
In this case, it induces an isomorphism
\[
 \Tr_\kk^G:\bar H_\kk(G)
       \xrightarrow{\;\sim\;}
       \cR_\kk(G)^*_{\mathrm{good}}.
\]
\end{theorem}

\begin{proof}
Suppose that $\ell>0$ is not banal for $G$.
Choose a compact open subgroup $K$ and an open
normal pro-$p$ subgroup $U$ such that
$\ell\mid[K:U]$.
Then $f=\mathbf 1_K$ satisfies
\[
 f*f=\mu_G(K)f=0,
\]
so its trace on every irreducible smooth
representation is zero.
On the other hand, evaluation at the identity
annihilates commutators, since $G$ is unimodular,
and $f(1)=1$.
Thus $f$ has a nonzero image in the cocenter,
and the trace map is not injective.
This proves the ``only if'' direction.

Suppose now that $\ell=0$ or that $\ell$ is banal
for $G$.
Under either assumption, every finite quotient of
a compact subgroup of $G$ has order invertible
in $\kk$.
This property passes to every Levi subgroup,
since its compact subgroups are contained in
compact open subgroups of $G$.

Fix a standard Levi subgroup $M$ and a split
maximal torus $T\subset M$.
The finite reductive quotient of a hyperspecial
subgroup of $M$ is split, and the normalizer of
a split maximal torus in this quotient surjects
onto $W_M$.
Hence $|W_M|\in\kk^\times$.

Put $Y=X_*(T)$, and let $\Omega_M$ be the
length-zero subgroup of $Y\rtimes W_M$.
As in \eqref{eq:inj-affine-central-index},
the natural projection to $W_M$ has kernel
$Y^{W_M}$ on $\Omega_M$, giving
\begin{equation}\label{eq:inj-split-central-index}
 \Omega_M/Y^{W_M}\hookrightarrow W_M,
 \qquad |\Omega_M/Y^{W_M}|\in\kk^\times.
\end{equation}
Choose a uniformizer $\varpi$ and set
\[
 Z_M=\{y(\varpi):y\in Y^{W_M}\}\subset Z(M),
 \qquad Q_M=M/Z_M.
\]
For every standard parahoric subgroup $\mathcal P$
of $M$, the group $N_M(\mathcal P)/Z_M$ is compact.
Its quotient by the image of $N_M(\mathcal P)^c$
has order dividing $|\Omega_M/Y^{W_M}|$.
Since finite quotients of $N_M(\mathcal P)^c$
also have orders invertible in $\kk$, the same
holds for finite quotients of
$N_M(\mathcal P)/Z_M$.
Consequently, all finite-dimensional Hecke algebras
occurring in the compact kernel formula for $Q_M$
are semisimple.
The formula in \S\ref{subsec:whole-kernel-proof}
therefore gives
\begin{equation}\label{eq:inj-split-quotient}
 \ker\bigl(\Tr_\kk^{Q_M}|_{\bar H_\kk(Q_M)^c}\bigr)
   =N_\kk(Q_M)=0.
\end{equation}

We now use the same argument with central
translations as in the proof of
\Cref{thm:inj-affine-whole}.
The pushforward
$\bar q_M:\bar H_\kk(M)\to\bar H_\kk(Q_M)$
identifies its target with the quotient by the
relations $z\cdot h-h$, for $z\in Z_M$.
Central translations shift Newton vectors by
their valuations, and the valuation homomorphism
is injective on $Z_M$.
Thus these translations act freely on the Newton
degrees, and, as in
\eqref{eq:inj-affine-degreewise}, $\bar q_M$
is injective on each Newton component.
Together with \eqref{eq:inj-split-quotient},
this makes every summand in the formula of
\Cref{thm:whole-kernel} zero.
Hence $\Tr_\kk^G$ is injective.

Finally, taking $M=G$ gives $|W_0|\in\kk^\times$.
The trace Paley--Wiener theorem
\Cref{thm:known-tpw} therefore gives surjectivity
onto the good forms.
\end{proof}

The same banality criterion also holds for certain
non-split groups.
For example, let $D/F$ be a central division algebra
of degree $d$, and let $q$ be the residue-field
cardinality of $F$.
For $\ell>0$, the trace map on the whole cocenter of $D^\times$ is injective if and only if $\ell$ is banal,
equivalently, $\ell\nmid(q^d-1)$.
Under banality, the Hecke algebras at the congruence
levels $1+\mathfrak p_D^r$ are finite direct sums
of matrix algebras over Laurent polynomial rings,
whose cocenters are separated by finite-dimensional
characters.

This example suggests that banality may characterize
injectivity on the whole cocenter for arbitrary
connected reductive $p$-adic groups.
We will explore this question in future work.

\section{Examples}

\label{sec:examples}

We illustrate the formulas for the kernel of the trace map for a split torus and for affine Hecke algebras of types $\widetilde A_2$ and $\widetilde C_2$. We use the notation of \S\ref{sec:preliminaries} for the torus and of \S\ref{sec:affine} for the affine examples.

\subsection{A split torus}
\label{subsec:torus-example}

In this example, assume $\ell\ne p$.

This example illustrates how the kernel on the
compact part extends to every Newton component.
Let $G=F^\times$, let $q$ be the residue-field
cardinality, and put $\CI_1=1+\mathfrak p_F$.
The quotient
$\mathcal O_F^\times/\CI_1\cong\mathbb F_q^\times$
is cyclic of order $q-1$.
Choosing a uniformizer and a generator of this
quotient gives
\[
 H_{1,\kk}\cong\kk[t,t^{-1}]\otimes_\kk B_\kk,
 \qquad B_\kk=\kk[c]/(c^{q-1}-1).
\]
Since $G$ is abelian, $H_{1,\kk}$ embeds in
$\bar H_\kk(G)=H_\kk(G)$.
Its Newton components are $t^rB_\kk$, $r\in\ZZ$,
and its compact part is $B_\kk$.

Suppose that $\ell>0$, and write $q-1=\ell^a d$
with $\ell\nmid d$.
Characters trivial on $\CI_1$ correspond to
evaluations at $t=u\in\kk^\times$ and $c=\zeta$
with $\zeta^d=1$.
Since $c^{q-1}-1=(c^d-1)^{\ell^a}$, we obtain
\[
 \ker(\Tr_\kk^G|_{H_{1,\kk}})
   =(c^d-1)H_{1,\kk}
   =\bigoplus_{r\in\ZZ}t^r(c^d-1)B_\kk.
\]
On each Newton component, the kernel is therefore
a translate of its kernel on the compact part
and has dimension $q-1-d$.

The subgroup $\langle c^{\ell^a}\rangle$ has order
$d$, and its group algebra maps isomorphically
onto $B_\kk/(c^d-1)B_\kk$.
Thus
\[
 B_\kk=(c^d-1)B_\kk
       \oplus\kk[\langle c^{\ell^a}\rangle].
\]
The primitive idempotents of the second summand
decompose $H_{1,\kk}$ into connected factors,
on each of which projective rank is constant.
Each idempotent is itself the rank of the
corresponding projective summand, so
\[
 \im\Rk_{H_{1,\kk}}
   =\kk[\langle c^{\ell^a}\rangle].
\]
This identifies the complementary subspace with
the image of the $\kk$-linearized rank map.
In characteristic zero, the kernel of the trace
map is zero and the image of the rank map is $B_\kk$.

\subsection{Affine type \texorpdfstring{$\widetilde A_2$}{A2}}
\label{subsec:affine-rank-two}

This example illustrates how the root datum and
the coefficient characteristic affect the kernel
of the trace map.
Consider the equal-parameter affine Hecke algebras
associated with $\mathrm{SL}_3$ and $\mathrm{PGL}_3$,
with parameter $q\in\kk^\times$.
Their Iwahori--Weyl groups are respectively
$Q^\vee\rtimes S_3$ and $P^\vee\rtimes S_3$,
where $Q^\vee\subseteq P^\vee$ are the coroot
and coweight lattices of type $A_2$.
Thus $\Omega=1$ for $\mathrm{SL}_3$, whereas
$\Omega=\langle\tau\rangle\cong C_3$ for
$\mathrm{PGL}_3$, with $\tau$ cyclically permuting
the affine simple reflections.

Set $a=\bar T_{s_i}$, independent of $i$, and
$c_{ij}=\overline{T_{s_i}T_{s_j}}$ for $i<j$.
By \Cref{thm:affine-cocenter-basis}, the rigid
cocenters have bases
\begin{equation}\label{eq:a2-rigid-basis}
 \begin{array}{ll}
 \mathrm{SL}_3:&1,\ a,\ c_{01},\ c_{12},\ c_{02},\\
 \mathrm{PGL}_3:&1,\ a,\ c,\ \overline{T_\tau},\
                         \overline{T_{\tau^2}},
 \end{array}
\end{equation}
for every characteristic and every nonzero parameter.
Here $c_{01}=c_{12}=c_{02}=c$ in the second case;
see also \cite[Examples~3.14--3.15]{He16}.

The three finite parahoric algebras indexed by
$\{0,1\}$, $\{1,2\}$, and $\{0,2\}$ are of type $A_2$.
Their character tables
\cite[Example~3.5]{He16} give the following
contributions to the kernel in every characteristic.
If $(q+1)(q^2+q+1)\ne0$, all three trace maps
are injective.
If $q=-1$, the images of their kernels in the
rigid cocenter coincide with $\kk(a+1)$.
If $q^2+q+1=0$, these images are generated by
\[
 z_{ij}=(q+1)c_{ij}+(q+2)a+1.
\]
Since $q+1\ne0$ in this case, their $c_{ij}$
coordinates show that the three classes are
linearly independent for $\mathrm{SL}_3$.
For $\mathrm{PGL}_3$, they coincide.

For $\mathrm{PGL}_3$, the additional finite
extended parahoric algebra $\kk[\Omega]$ contributes
nothing when $\ell\ne3$.
When $\ell=3$, the image of the kernel of its
trace map is
\[
 \spanop_\kk\{
   \overline{T_\tau}-1,\,
   \overline{T_{\tau^2}}-1
 \}.
\]
By \eqref{eq:a2-rigid-basis}, these two classes
are independent modulo the contributions from
the type $A_2$ algebras.
They therefore increase the dimension of the
kernel by two for every $q$.
These relations also show that the kernel need
not respect the $C_3$-grading in characteristic $3$.

At $q=1$, the formulas recover
\Cref{lem:finite-regular}:
$a+1=a-1$ in characteristic $2$, and
$z_{ij}=-(c_{ij}-1)$ in characteristic $3$.

By \Cref{thm:affine-kernel}, the contributions above span
the kernel of the trace map on each rigid cocenter.
Since both rigid cocenters have dimension five, subtracting
the dimension of each kernel gives the ranks in the table.

\begin{center}
\begin{tabular}{@{}lccc@{}}
\toprule
Parameter
 & $\mathrm{SL}_3$
 & \multicolumn{2}{c}{$\mathrm{PGL}_3$}\\
\cmidrule(l){3-4}
 & all $\ell$ & $\ell\ne3$ & $\ell=3$\\
\midrule
$(q+1)(q^2+q+1)\ne0$ & $5$ & $5$ & $3$\\
$q=-1$               & $4$ & $4$ & $2$\\
$q^2+q+1=0$          & $2$ & $4$ & $2$\\
\bottomrule
\end{tabular}
\end{center}

When $\ell\notin\{2,3\}$, the order of $W_0$ is invertible
in $\kk$. By \Cref{cor:affine-character-extension}, each rank
in the table then equals the dimension of the corresponding
rigid quotient. In particular, this confirms the dimension
predictions in \cite[Examples~3.14--3.15]{He16} over $\CC$.

\subsection{Unequal parameters in affine type
\texorpdfstring{$\widetilde C_2$}{C2}}
\label{subsec:affine-c2}

For unequal parameters, semisimplicity of the finite
Hecke algebra associated with $W_0$ does not suffice
for injectivity of the trace map.
Consider affine type $\widetilde C_2$ with $\Omega=1$
over $\CC$, and take
\[
 (q_0,q_1,q_2)=(-2,2,3).
\]
All rank-one parahoric algebras are semisimple.
By the finite character calculation in
\cite[Example~7.9]{CH}, the algebra associated with
$W_0=\langle s_1,s_2\rangle$ is also semisimple,
whereas the algebra on $\{0,1\}$ is not, since
\[
 (q_1+q_2)(1+q_1q_2)=35\ne0,
 \qquad
 (q_0+q_1)(1+q_0q_1)=0.
\]
Thus \Cref{cor:affine-injectivity} shows that
the trace map is not injective on the compact
cocenter.

\end{document}